%% file: main.tex
\documentclass[a4paper,12pt]{amsart}

\input{preamble}

\title{Limits in categories of \'etale groupoids and pseudogroups}

\author{Jonathan Taylor}
\address{Jonathan Taylor, Institut für Mathematik, University of Potsdam, Campus Golm, Haus 9, Karl-Liebknecht-Str. 24--25, 14476, Germany}
\email{jonathan.taylor@uni-potsdam.de}

\begin{document}
	\begin{abstract}
		We show that the category of sober \'etale groupoids and actors admits all small limits.
		This is achieved by computing the limits in the equivalent category of pseudogroups with pseudogroup morphisms, which we show admits a forgetful functor to the category of sets which creates limits.
		We give an alternative proof of the adjunction of Cockett and Garner \cite{CockettGarner_GeneralisingEtaleGrpds} in the specific setting of \'etale groupoids and pseudogroups which is a central tool for computing limits of sober \'etale groupoids.
	\end{abstract}

	\maketitle
	\pagenumbering{arabic}
	
	\section{Introduction}
	\subfile{sections/intro}
	
	\section{Preliminaries}\label{sec-prelims}
	\subfile{sections/prelims}

	\section{Functors between groupoids and pseudogroups}\label{sec-cats}
	\subfile{sections/cats}

	\section{The adjunction}\label{sec-adjunction}
	\subfile{sections/adjunction}

	\section{Limits}\label{sec-lims}
	\subfile{sections/lims}

	\appendix

	\section*{Appendices}
	As the target audience of this article includes readers who may not specialise in category theory or Stone duality, we include appendices containing key definitions and results for each of these topics.
	\section{Category theory}\label{app-CatTheory}
	\subfile{appendices/catTheory}

	\section{Stone duality}\label{app-StoneDuality}
	\subfile{appendices/stoneDuality}

	\printbibliography
	
\end{document}

%% file: preamble.tex
\usepackage[margin=3cm]{geometry}
\usepackage{amssymb}
\usepackage{amsmath}
\usepackage{amsthm}
\usepackage{mathrsfs}
\usepackage{enumitem}
\usepackage{hyperref}

\usepackage[
    backend=biber,
    style=numeric
]{biblatex}
\usepackage{tikz}
\usepackage{tikz-cd}
\usetikzlibrary{matrix,calc,arrows,decorations.pathmorphing}
\tikzset{node distance=2cm, auto}

\usepackage{subfiles}

\tikzset{cd/.style=matrix of math nodes,row sep=2em,column sep=2em, text height=1.5ex, text depth=0.5ex}
\tikzset{cdar/.style=->,auto}
\tikzset{mid/.style={anchor=mid}} 
\tikzset{narrowfill/.style={inner sep=1pt, fill=white}}

\newtheorem{theorem}{Theorem}[section]
\newtheorem{lemma}[theorem]{Lemma}
\newtheorem{corollary}[theorem]{Corollary}

\newtheorem{proposition}[theorem]{Proposition}
\newtheorem*{theorem*}{Theorem}

\theoremstyle{definition}
\newtheorem{definition}[theorem]{Definition}
\newtheorem{example}[theorem]{Example}
\newtheorem{notation}[theorem]{Notation}

\theoremstyle{remark}
\newtheorem{remark}[theorem]{Remark}

\newcommand{\NN}{\mathbb{N}}

\newcommand{\ZZ}{\mathbb{Z}}

\newcommand{\Cc}{\mathcal{C}}
\newcommand{\Dd}{\mathcal{D}}
\newcommand{\Ee}{\mathcal{E}}

\newcommand{\Oo}{\mathcal{O}}
\newcommand{\Pp}{\mathcal{P}}

\newcommand{\Uu}{\mathcal{U}}

\newcommand{\Set}{\mathbf{Set}}
\newcommand{\Top}{\mathbf{Top}}

\newcommand{\Grp}{\mathbf{Grp}}

\newcommand{\Ring}{\mathbf{Ring}}

\newcommand{\Vect}{\mathbf{Vect}}
\newcommand{\Frame}{\mathbf{Frame}}

\newcommand{\EGA}{\mathbf{\acute{E}tAct}}

\newcommand{\Pseu}{\Psi\mathbf{Grp}}
\newcommand{\Nat}{\operatorname{Nat}}

\newcommand{\colim}{\operatorname{colim}}

\newcommand{\op}{\operatorname{op}}

\newcommand{\Id}{\operatorname{id}}

\newcommand{\Sob}{\mathrm{Sob}}
\newcommand{\Spat}{\mathrm{Spat}}
\renewcommand{\op}{{\mathrm{op}}}
\newcommand{\Bis}{\operatorname{Bis}}
\newcommand{\dom}{\operatorname{dom}}
\newcommand{\Cone}{\operatorname{Cone}}
\newcommand{\coeq}{\operatorname{coeq}}

\newcommand{\baltimes}[2]{
	\mathbin{_{#1}\times_{#2}}%
}

%% file: sections/intro.tex
An important generalisation of topological spaces is that of topological groupoids which combine geometric, algebraic, and dynamical structure into one unifying class of objects containing groups, spaces, equivalence relations, and many other objects of interest in many areas of mathematics.
Topological groupoids have seen successful application in operator algebras.
Many classes of interesting \(C^*\)-algebras arise as groupoid \(C^*\)-algebras, including graph \(C^*\)-algebras, group \(C^*\)-algebras, crossed products of commutative \(C^*\)-algebras.
Of particular interest are \'etale groupoids: groupoids for which the source and range structure maps are local homeomorphisms.
The fibrewise-discrete structure of \'etale groupoids allows techniques common in the study of discrete (or simply untopologised) groups to be applied and generalised, while still retaining great use in operator algebras.

The construction of a \(C^*\)-algebra from a groupoid is however not functorial for groupoid homomorphisms (i.e. continuous functors).
The philosophical obstruction to such a functoriality is that groupoid \(C^*\)-algebras generalise both group \(C^*\)-algebras and \(C^*\)-algebras associated to (locally compact Hausdorff) topological spaces.
The group \(C^*\)-algebra construction is covariantly functorial for group homomorphisms, whereas the corresponding construction for topological spaces is contravariantly functorial for continuous maps by Gelfand duality.

To resolve this, one may consider a different class of morphisms between groupoids.
Buneci and Stachura \cite{BuneciStachura_morphismslocallycompactgroupoids} showed that \emph{actors} between groupoids lift covariantly to groupoid \(C^*\)-algebras generalising both the universality for group \(C^*\)-algebras and Gelfand duality.
The resulting category of groupoids with actors then lends itself naturally to the study of groupoid \(C^*\)-algebras, and moreover all invertible actors arise from isomorphisms of groupoids as topological categories.
A natural question from the operator algebraic perspective is to then ask how well-behaved this functor from groupoids to \(C^*\)-algebras is.
In particular, a motivating question is whether this functor preserves limits, as this would allow one to construct more groupoid \(C^*\)-algebras with universal properties from groupoids with universal properties.
To this end, this article seeks to study limits in the category of groupoids and actors.

The primary discomfort that comes with the category of groupoids and actors is that it is not concrete (at least, not obviously): there is no obvious faithful forgetful functor to the category of sets.
Constructions of limits and colimits in this category become opaque in turn, as we lose access to the standard strategy of describing morphisms as functions between underlying sets.
To this end, we consider pseudogroups as an ersatz for groupoids and actors which retain enough information to reconstruct groupoids and actors while being more amenable to concrete constructions.
\'Etale groupoids admit a base of open subsets on which the range and source maps are homeomorphisms onto their images.
Such open subsets are called \emph{bisections}, and instead of studying a groupoid from its points, we may study it from its bisections.
The collection of bisections of an \'etale groupoid forms a pseudogroup: a complete inverse semigroup with respect to a natural order structure.
Actors between groupoids then induce homomorphisms between their pseudogroups of bisections, and this construction is functorial.

In order to study \'etale groupoids from their pseudogroups of bisections, we require a duality between \'etale groupoids and pseudogroups.
To begin, we need a method of reconstructing a topological space from the algebraic structure of its collection of open sets.
The theory of frames provides a point-agnostic approach to topology, focussing on the lattice structure of open sets rather than an underlying set of points.
Of course not every topological space may be recovered from its frame of open sets, but a surprisingly large class of spaces may be.
Such spaces are called \emph{sober}, and the class of sober spaces is large enough to include all Hausdorff spaces, as well as any space admitting a local homeomorphism to a sober space.
The fundamental result in frame theory linking it to topology is the duality between sober spaces and spatial frames (frames which are isomorphic to topologies as lattices): there is an adjunction between the category of topological spaces and the category of frames which restricts to an equivalence between the full subcategories of sober spaces and spatial frames.

The duality between sober spaces and spatial frames was extended to a duality between \'etale groupoids and pseudogroups by Matsnev and Resende \cite{MatsnevResende_EtaleGrpdsGermGrpds}, showing that sober \'etale groupoids can be reconstructed from their pseudogroups of bisections.
Lawson and Lenz \cite{LawsonLenz_PseuGrpsEtaleGrpds} extended this duality to some classes of morphisms, considering \emph{callitic} pseudogroup homomorphisms and \emph{continuous covering functors} between \'etale groupoids.
Going even further, Cockett and Garner \cite{CockettGarner_GeneralisingEtaleGrpds} showed that the category of \'etale partite internal groupoids to a join restriction category (a generalisation of \'etale groupoids) and partite internal cofunctors (a generalisation of actors) is equivalent to the category of join restriction categories and hyperconnected join restriction functors.
In particular, this last generalisation allows us to safely replace sober \'etale groupoids and actors between them by their pseudogroups of bisections and pseudogroup homomorphisms.

The purpose of this article is to utilise the adjunction between \'etale groupoids and pseudogroups to compute limits in the category of \'etale groupoids and actors.
This is achieved by showing that the category of sober \'etale groupoids with actors is equivalent to the category of spatial pseudogroups with pseudogroup homomorphisms.

\begin{theorem*}[{Theorem~\ref{thm-adjunction}, Corollary~\ref{cor-catEquivSobSpat}}]
    The spatialisation functor \(\Sigma\) sending a pseudogroup to its spatialisation groupoid is left-adjoint to the bisection functor \(\Bis\) sending an \'etale groupoid to its pseudogroup of bisections.
    This adjunction restricts to an equivalence of categories between sober \'etale groupoids and spatial pseudogroups.
\end{theorem*}

This adjunction and the resulting equivalence of categories is not a new result: the version presented in this article can be attained as an application of \cite[Theorem~6.8]{CockettGarner_GeneralisingEtaleGrpds}.
We present a proof of the theorem which is not dependent of the sophisticated machinery developed by Cockett and Garner in \cite{CockettGarner_GeneralisingEtaleGrpds} in the interest of making the result more accessible to non-experts.

To compute limits, we show that forgetful functor from (spatial) pseudogroups to sets creates limits, so implies completeness of the category of pseudogroups.

\begin{theorem*}[{Theorems~\ref{thm-pseuComplete}, \ref{thm-spatSobComplete}}]
    The forgetful functor from the category of (spatial) pseudogroups to sets creates limits.
\end{theorem*}

The above two results then imply that the category of sober \'etale groupoids and actors is complete, and that limits may be computed by first computing the limit in the category of sets, then lifting to a pseudogroup, and finally taking the spatialisation groupoid of this pseudogroup.

The article is structured as follows.
Section~\ref{sec-prelims} provides the necessary preliminaries on groupoids, inverse semigroups, and pseudogroups.
Section~\ref{sec-cats} describes two functors which form the adjunction between \'etale groupoids and pseudogroups: the bisection functor sending a groupoid to its pseudogroup of bisections, and the spatialisation functor sending a pseudogroup to its associated \'etale groupoid.
Section~\ref{sec-adjunction} then shows that these two functors form an adjunction between the categories of \'etale groupoids and pseudogroups, generalising the adjunction between topological spaces and frames.
This adjunction is then shown to restrict to an equivalence of categories between sober groupoids and spatial pseudogroups.
Finally, in Section~\ref{sec-lims} we show that the forgetful functor from pseudogroups to sets creates limits, and hence that the category of pseudogroups is complete.
This then implies that the category of sober \'etale groupoids is complete, and we explicitly describe products and equalisers in this category.
We compute some examples of limits, including the following pullback graph groupoids.
We further investigate how the inclusion functors mapping groups to both groupoids and pseudogroups interact with limits.
We also briefly discuss colimits in this category, and some of the issues that may arise.

We have also included appendices with key definitions and results in category theory and in Stone duality.

%% file: sections/prelims.tex
The main tool of this article is a generalisation of the duality between (sober) topological spaces and (spatial) frames.
For a detailed exploration of the theory of frames and locales, we refer the reader to \cite{PicadoPultr_FramesLocales}.
We provide an abridged and mostly self-contained recollection of the main results necessary for the article in Appendix~\ref{app-StoneDuality}.
We also recall elementary and useful properties of groupoids, inverse semigroups, and pseudogroups, as well as set notation.

Our notation for the application of a functor \(F\) to an object \(x\) or a morphism \(f\) may be denoted without brackets e.g. \(Fx\) and \(Ff\).
For a topological space \(X\), the topology will be denoted by \(\Oo X\). 
A number of relevant elementary definitions and results in category theory are given in Appendix~\ref{app-CatTheory}.

\subsection{Topological spaces and frames}

The definition of a frame is an abstract axiomatisation of the lattice structure of open sets in a topological space.
We give some of the basic definitions and relevant results for frame theory and the relationship to topology here.

\begin{definition}
    A \emph{frame} \(F\) is a partially ordered set with finite infima and all suprema, in which infima distribute over suprema.
    We denote the infimum or \emph{meet} of two elements \(e,f\in F\) by \(e\wedge f\), and the supremum or \emph{join} of a family \((e_\alpha)_\alpha\subseteq F\) by \(\bigvee_\alpha e_\alpha\).
    Distributivity is then expressed as 
    \[e\wedge\bigvee_\alpha f_\alpha=\bigvee_\alpha (e\wedge f_\alpha)\]
    for all \(e,f_\alpha\in F\).

    Given frames \(F\) and \(L\), a \emph{frame homomorphism} \(\varphi\colon F\to L\) is a function which preserves finite infima and arbitrary suprema. 
    Frames together with frame homomorphisms form the category \(\Frame\).
\end{definition}

For any topological space \(X\), the topology \(\Oo X=\{U\subseteq X\text{ open}\}\) is an example of a frame. 
Frames which are isomorphic to \(\Oo X\) for some topological space \(X\) are called \emph{spatial}.
Moreover, continuous maps \(f\colon X\to Y\) induce frame homomorphisms \(f^{-1}\colon\Oo Y\to\Oo X\), describing a functor \(\Oo\colon \Top^\op\to\Frame\), where \(\Top\) is the category of topological spaces.

\begin{definition}\label{defn-frameChar}
    A \emph{character} on a frame \(F\) is a frame homomorphism \(\chi\colon F\to\Omega\).
    Denote the set of characters by \(\widehat F\).
\end{definition}

We topologise \(\widehat F\) by open sets \(\Uu_{e}=\{\chi\in\widehat F: \chi(e)=1\}\) for \(e\in F\).
Any frame homomorphism \(L\to F\) induces a continuous map \(\widehat F\to\widehat L\) via precomposition.
This assignment describes a functor \(\sigma\colon\Frame\to\Top^\op\) with \(\sigma F=\widehat F\) for a frame \(F\).

A topological space is \emph{sober} if every closed irreducible subset is the closure of a unique singleton (see Definition~\ref{defn-irredSob} for more details).
Stone duality is the assertion that sober spaces and spatial frames form equivalent categories.
We denote by \(\Frame_\Spat\) and \(\Top_\Sob\) the full subcategories of \(\Frame\) and \(\Top\) spanned respectively by spatial frames and sober spaces.

\begin{theorem}[see Theorem~\ref{thm-appSoberSpatAdj}]\label{thm-soberSpatAdj}
    There is an adjunction \(\sigma\dashv\Oo\) which restricts to an equivalence of categories \(\Frame_\Spat\simeq\Top_\Sob^\op\).
\end{theorem}

We will later be concerned with \emph{\'etale groupoids}, which come equipped with local homeomorphisms to an open subspace (the \emph{unit space} of the groupoid).
The next proposition shows that spaces which are locally homeomorphic to sober spaces are themselves sober, so sobriety can in this way be verified locally.

\begin{proposition}\label{prop-sobLocHomeo}
    Let \(X\) and \(Y\) be topological spaces and suppose \(s\colon X\to Y\) is a local homeomorphism.
    If \(Y\) is sober, then so is \(X\).
    \begin{proof}
        Let \(C\subseteq X\) be a closed irreducible subset.
        Since \(s\) is a local homeomorphism, there is an open subset \(U\subseteq X\) with \(U\cap C\neq\emptyset\) on which \(s\) is injective, and hence is a homeomorphism onto its open image in \(Y\).
        The subset \(U\) is sober since it is homeomorphic to an open subset of a sober space.
        Moreover, since \(C\cap U\) is an irreducible subset of \(U\), its image \(s(C\cap U)\) is irreducible in \(s(U)\) and so there is a unique point \(y\in s(U)\) with \(\overline{\{y\}}=s(C\cap U)\).

        The set \(\overline{s(C)}\) is closed and irreducible in \(Y\) since any two closed subsets \(D_1,D_2\subseteq Y\) with \(D_1\cup D_2=\overline{s(C)}\) lift to closed subsets \(C_i:=s^{-1}(D_i)\cap C\) with \(C_1\cup C_2=C\), whereby \(C=C_i\) for some \(i\in\{1,2\}\).
        Then \(s(C)=s(C_i)\subseteq s(s^{-1}(D_i))\subseteq D_i\subseteq \overline{s(C)}\), so \(\overline{s(C)}=D_i\) as \(D_i\) is closed.Thus \(\overline{s(C)}\) is irreducible, and sobriety of \(Y\) yields a point \(y\in Y\) with \(\overline{\{y\}}=\overline{s(C)}\).
        
        We claim there is a point \(x\in C\) with \(s(x)=y\).
        Suppose this is not the case.
        For each \(z\in X\) with \(s(z)=y\) we may find an open neighbourhood \(U_z\subseteq X\) of \(z\) onto which \(s\) restricts to an injective map.
        Since \(C\) is closed and does not contain \(z\) we may further shrink \(U_z\) to ensure \(U_z\cap C=\emptyset\) for all such \(z\).
        Then \(V:=\bigcup_{z\in s^{-1}(\{y\})}U_z\) is an open subset of \(X\) containing the fibre \(s^{-1}(\{y\})\), and 
        \[s(V)\cap s(C)=\bigcup_{z\in s^{-1}(\{y\})}s(U_z)\cap s(C)=\bigcup_{z\in s^{-1}(\{y\})}s(U_z\cap C)=\emptyset,\]
        since \(s\) is injective on each \(U_z\).
        But then \(s(V)\) is an open neighbourhood of \(y\) in \(Y\) which does not intersect \(s(C)\), contradicting that \(s(C)\) is dense in \(\overline{s(C)}=\overline{\{y\}}\).
        Thus we may select \(x\in C\) with \(s(x)=y\).

        Let \(U\subseteq X\) be an open neighbourhood of \(x\) onto which \(s\) restricts to an injective map.
        The set \(C\cap U\) is an irreducible closed subset of \(U\) in the subspace topology.
        Indeed, closed subsets of \(U\) have the form \(D\cap U\) for some closed \(D\subseteq X\), and if \((C_1\cap U)\cup(C_2\cap U)=C\cap U\) for some \(C_1,C_2\subseteq X\) closed, then 
        \[C=(C_1\cap C)\cup (C_2\cap C)\cup (X\setminus U\cap C),\]
        whereby the irreducibility of \(C\) implies either \(C_1\cap C\), \(C_2\cap C\), or \(X\setminus U\cap C\) is equal to \(C\).
        The last case can be excluded since \(X\setminus U\cap C\) cannot be all of \(C\) since \(U\cap C\neq \emptyset\).
        Hence \(C_i\cap C=C\) for some \(i=1,2\) and we have \(C\cap U=C_i\cap U\), hence \(C\cap U\) is irreducible in \(U\).
        Since \(U\cap C\) is an open subset of \(C\) (again in the subspace topology), it is necessarily dense (otherwise \(C\setminus U\cup \overline{C\cap U}=C\) is a non-trivial reduction).
        
        Since \(C\cap U\) is dense in \(C\) we have
        \[C=\overline{C\cap U}=\overline{\{x\}}\subseteq C\]
        noting that \(x\in C\) implies the closure of the singleton \(\{x\}\) is contained in \(C\).
        
        It remains to show that the point \(x\) is unique.
        To this end, it suffices to show that \(X\) is \(T_0\), since any two distinct singletons in a \(T_0\) space have distinct closures.
        Fix \(x_1,x_2\in X\) and an open neighbourhood \(U\subseteq X\) of \(x_1\) onto which \(s\) restricts to a homeomorphism onto its image.
        If \(x_2\notin U\) then we are done.
        If \(s(x_1)=s(x_2)\) either \(x_1=x_2\) or \(x_2\notin U\) since \(s\) is injective on \(U\).
        Lastly, if \(s(x_1)\neq s(x_2)\) then (since \(Y\) is \(T_0\)) we may find an open subset \(V\subseteq Y\) with either \(s(x_1)\in V\) and \(s(x_2)\notin V\) or vice versa.
        Then either \(s^{-1}(U\cap V)\) contains \(x_1\) but not \(x_2\), or \(s^{-1}(V)\) contains \(x_2\) but not \(x_1\), yielding \(T_0\).
    \end{proof}
\end{proposition}

\subsection{\'Etale groupoids and actors}

By \emph{groupoid} we always mean a small category in which all morphisms are invertible.
We denote groupoids by \(G,H,K\) and morphisms of groupoids are referred to as \emph{arrows} and typically denoted by \(\gamma,\eta,x,y\). 
For a groupoid \(G\), we shall consider the unit space \(G^{(0)}\) as a subset of the set of arrows (which we denote by \(G\)) by identifying objects with their corresponding identity morphisms.
The set of composable pairs in \(G\) is denoted by \(G^{(2)}\), and the composition map \(G^{(2)}\to G\) is denoted by juxtaposition \((\gamma,\eta)\mapsto\gamma\eta\). 
The inversion map in \(G\) is denoted by \(\gamma\mapsto\gamma^{-1}\).
The range and source maps \(r,s\colon G\to G^{(0)}\) send an arrow \(\gamma\) each to \(r(\gamma)=\gamma\gamma^{-1}\) and \(s(\gamma)=\gamma^{-1}\gamma\).
In particular, the composition of arrows is read right-to-left, following the commonly used convention for function composition.

A \emph{topological groupoid} is a groupoid \(G\) equipped with a topology \(\Oo G\) for which the inverse operation and multiplication are all continuous (where \(G^{(2)}\) is equipped with the subspace topology inherited from the product topology of \(G\times G\)).
A topological groupoid \(G\) is \emph{\'etale} if the range and source maps are local homeomorphisms (note that either of these maps being a local homeomorphism implies the other is as well).
We shall often invoke the equivalent characterisation that \(G\) is \'etale if and only if the topology of \(G\) admits a base of \emph{bisections}; open subsets \(U\subseteq G\) to which the range and source maps restrict to homeomorphisms onto their images. 

The ability to restrict to open bisections in \'etale groupoids is invaluable in their study, and much of the analysis is achieved by leveraging the fact that an \'etale groupoid is locally homeomorphic to its unit space. 
Proposition~\ref{prop-sobLocHomeo} demonstrates it suffices to show that the unit space of an \'etale groupoid is sober to conclude that the entire space is sober.

We now introduce the morphisms we consider between \'etale groupoids.
There is a modest zoo of names for these morphisms between groupoids, including Zakrzewski morphisms \cite{BuneciStachura_morphismslocallycompactgroupoids}, comorphisms \cite{HigginsMackenzie_DualityMorphismVectBunds}, algebraic morphisms \cite{BuneciStachura_morphismslocallycompactgroupoids}, actors \cite{MeyerZhu_Groupoids}, and cofunctors \cite{CockettGarner_GeneralisingEtaleGrpds}.
We adopt the terminology of \emph{actor} to emphasise the similarity and relationship to groupoid actions.

Setting some notation: for (continuous) functions \(f\colon X\to Z\) and \(g\colon Y\to Z\) we denote the pullback of \(f\) and \(g\) by \(X\baltimes{f}{g}Y:=\{(x,y)\in X\times Y: f(x)=g(y)\}\).
When applicable, we equip the pullback with the subspace topology inherited from the product topology of \(X\times Y\).

\begin{definition}[{\cite[23]{MeyerZhu_Groupoids}, \cite[Definition~2.2]{Tay_FunctGrpdCstarAlg}}]\label{defn-actor}
    Let \(G\) and \(H\) be \'etale groupoids.
    An \emph{actor} \(h\colon G\curvearrowright H\) consists of a continuous map \(\rho_h\colon H\to G^{(0)}\) and a continuous map \(\cdot_h: G\baltimes{s}{\rho}H\to H\), \((g,x)\mapsto g\cdot_h x\), satisfying the following conditions:
    \begin{enumerate}
        \item \(\rho_h(x)\cdot_h x=x\) for all \(x\in H\);
        \item \(r(\gamma)=\rho_h(\gamma\cdot_h x)\) for all \((\gamma,x)\in G\baltimes{s}{\rho}H\);
        \item \(\gamma_1\cdot_h(\gamma_2\cdot_h x)=(\gamma_1\gamma_2)\cdot_hx\) for all \((\gamma_1,\gamma_2)\in G^{(2)}\), \((\gamma,x)\in G\baltimes{s}{\rho} H\);
        \item \(\rho_h(x_1x_2)=\rho_h(x_1)\) for all \((x_1,x_2)\in H^{(2)}\);
        \item \(\gamma\cdot_h(x_1x_2)=(\gamma\cdot_h x_1)x_2\) for all \((\gamma,x_1)\in G\baltimes{s}{\rho}H\), \((x_1,x_2)\in H^{(2)}\).
    \end{enumerate}
    The map \(\rho_h\) is called the \emph{anchor}, and the function \(\cdot_h\) is called the \emph{actor multiplication}.
    We may omit the subscript \(h\) when it is clear from context.
\end{definition}

The first three conditions in Definition~\ref{defn-actor} describe a left-action of \(G\) on \(H\), treating \(H\) as a topological space.
The final two conditions state that the left action should commute with the right action of \(H\) on itself by right-multiplication.

The other typical morphisms considered between (topological) groupoids are (continuous) functors.
Meyer and Zhu show that the two concepts align for continuous functors which restrict to homeomorphisms of the unit spaces, respectively for actors whose anchor maps are homeomorphisms (see \cite[Example~4.17]{MeyerZhu_Groupoids}).

Given actors \(h\colon G\curvearrowright H\) and \(k\colon H\curvearrowright K\) with respective anchor maps \(\rho_h\) and \(\rho_k\), the composition \(kh\) is defined with anchor map \(\rho_{kh}:=\rho_k\circ\rho_h\) and multiplication \(\gamma\cdot_{kh}t:=(\gamma\cdot_h\rho_k(t))\cdot_k t\) for \(\gamma\in G\) and \(t\in K\) with \(s(\gamma)=\rho_{kh}(t)\).
One readily verifies that the composition \(kh\) is an actor from \(G\) to \(K\).
The left-multiplication of a groupoid \(G\) on itself (with range map as the anchor) acts as the identity actor at \(G\), forming a category.

\begin{definition}
    The category \(\EGA\) consists of \'etale groupoids as objects and actors as morphisms with the above described composition.
    We denote the full subcategory spanned by sober \'etale groupoids by \(\EGA_\Sob\).
\end{definition}

Isomorphisms in \(\EGA\) are in bijective correspondence with the `usual' topological groupoid isomorphisms; homeomorphic isomorphisms of the categories (see \cite[Proposition~4.19]{MeyerZhu_Groupoids}).
Thus two groupoids are isomorphic in \(\EGA\) if and only if they are isomorphic as topological categories, so we may safely refer to two `isomorphic' groupoids without specifying the ambient category.

The data of an actor \(X\curvearrowright Y\) for two topological spaces \(X\) and \(Y\) (viewed as groupoids consisting only of units) reduces to the anchor map \(\rho\colon Y\to X\), since units must act trivially by the first axiom in Definition~\ref{defn-actor}.
Conversely, every continuous function \(Y\to X\) is the anchor map for such an actor, so \(\Top^\op\) embeds into \(\EGA\) as a full subcategory.

\begin{notation}\label{not-Ugamma}
    Let \(G\) be an \'etale groupoid.
    For an element \(\gamma\in G\) and an open bisection \(U\in\Bis G\), if \(s(\gamma)\in r(U)\), then there is a unique element \(\gamma'\in U\) such that \((\gamma,\gamma')\) is a composable pair.
    We denote the composition of this pair by \(\gamma U\).
    Similarly, if \(r(\gamma)\in s(U)\), we denote by \(U\gamma\) the product \(\gamma'\gamma\), where \(\gamma'\) is the unique element of \(U\) with \(s(\gamma')=r(\gamma)\).
    More generally, if \(h\colon G\curvearrowright H\) is an actor, then for \(U\in\Bis G\) and \(x\in H\) with \(\rho(x)\in s(U)\), we denote by \(U\cdot_h x\) (or simply \(U\cdot x\) if the actor \(h\) is understood) the element \(\gamma\cdot x\), where \(\gamma\in U\) is the unique element of \(U\) with \(s(U)=\rho(x)\).
\end{notation}

An actor \(h\colon G\curvearrowright H\) induces a functor from the category of \(H\)-actions to the category of \(G\)-actions which leaves the underlying spaces invariant (see \cite[Proposition~4.18]{MeyerZhu_Groupoids}).
Continuous functors between groupoids do not have this property, which is one motivation for considering actors as morphisms instead of functors.
Another motivation comes from operator algebras: an actor \(h\colon G\curvearrowright H\) between \'etale groupoids induces a \({}^*\)-homomorphism \(C^*(G)\to C^*(H)\) between the \(C^*\)-algebras of the groupoids (see \cite{BuneciStachura_morphismslocallycompactgroupoids}, \cite{Tay_FunctGrpdCstarAlg}), yielding a functor from \(\EGA\) into the category of \(C^*\)-algebras.
While this functor is far from being full in general, the \({}^*\)-homomorphisms which arise from this construction encode dynamical information between the \(C^*\)-algebras, and in some circumstances an actor may be (re)constructed from a \({}^*\)-homomorphism (see \cite[Proposition~5.7]{Tay_FunctGrpdCstarAlg}).

\subsection{Pseudogroups}

Viewing groupoids as a generalisation of topological spaces by way of adding algebraic or dynamical information, the lens of Stone duality then raises the question: what is the corresponding generalisation of frames to encode this new information?
This question already has a number of answers in various forms (including but not limited to \cite{Resende_EtaleGrpdsQuantales}, \cite{LawsonLenz_PseuGrpsEtaleGrpds}, \cite{BussExelMeyer_InvSemiGrpActionsGrpdActions}, \cite{CockettGarner_GeneralisingEtaleGrpds}).

We introduce inverse semigroups and pseudogroups following \cite{Resende_EtaleGrpdsQuantales}, \cite{LawsonLenz_PseuGrpsEtaleGrpds} and will prove a version of Stone duality for \'etale groupoids in Section~\ref{sec-adjunction}.

\begin{definition}\label{defn-invSemGrp}
    An \emph{inverse semigroup} \(S\) is a semigroup with the following property: for each \(t\in S\) there is a unique element \(t^*\in S\) satisfying 
    \[t=tt^*t\text{, and }t^*=t^*tt^*.\]
\end{definition}

We refer the reader to \cite{Lawson_InvSemiGrps} for a general reference on the theory of inverse semigroups.
There is a natural partial ordering on an inverse semigroup given by defining \(u\leq t\) if \(u=tu^*u\).
The subset of idempotent elements of \(S\) is denoted by \(ES\), and this forms a meet-semilattice with respect to this ordering, where the meet operation is given by multiplication. 
The uniqueness criterion for the element \(t^*\) respective to \(t\) ensures that \(ES\) is a commutative subsemigroup of \(S\) (in fact, this is equivalent to the uniqueness of the element \(t^*\) for each \(t\in S\)).
This uniqueness also implies that any semigroup homomorphism between inverse semigroups is automatically \({}^*\)-preserving.

We briefly recall some of the elementary properties of the natural partial order on inverse semigroups which will be useful later in the article.

\begin{lemma}\label{lem-starOrdPres}
    The \({}^*\)-operation on an inverse semigroup \(S\) is an order automorphism.
    \begin{proof}
        The \({}^*\)-operation is an involution since the partial inverses are unique, so we need only show that it is order preserving.
        Fix \(t\leq u\in S\).
        Then \(t^*=(tu^*u)^*=u^*ut^*\), whereby \(t^*\leq u^*\) by \cite[1.4.6]{Lawson_InvSemiGrps}.
    \end{proof}
\end{lemma}

\begin{lemma}\label{lem-invSGrpHomOrdPres}
    Let \(\varphi\colon S\to T\) be a semigroup homomorphism between inverse semigroups.
    Then \(\varphi\) is \({}^*\)-preserving and order-preserving.
    \begin{proof}
        The uniqueness of partial inverses in \(S\) and \(T\) ensure that \(\varphi\) is \({}^*\)-preserving.
        For \(t\leq u\) we have \(\varphi(t)=\varphi(ut^*t)=\varphi(u)\varphi(t)^*\varphi(t)\), whereby \(\varphi(t)\leq\varphi(u)\).
    \end{proof}
\end{lemma}

Although inverse semigroups carry an order structure and an involution, this is all encoded in the semigroup operation as Lemma~\ref{lem-invSGrpHomOrdPres} shows.
Hence all semigroup homomorphisms automatically preserve these structures, so we need not add hypotheses to our morphisms to account for this.

Our main motivating example of an inverse semigroup comes from the theory of \'etale groupoids.

\begin{example}\label{ex-bisGInvSemGrp}
    Let \(G\) be an \'etale groupoid.
    The set 
    \[\Bis G:=\{U\subseteq G: U\text{ is an open bisection}\}\] 
    forms a semigroup with respect to the operation 
    \[UV:=\{\gamma\eta: \gamma\in U,\eta\in V, s(\gamma)=r(\eta)\}.\]
    It is easy to see that the range and source maps in \(G\) are injective on \(UV\), since they each are on \(U\) and \(V\) individually.
    Since the range map in \(G\) is open, the multiplication map is as well by \cite[Lemma~2.4.11]{Sims_EtaleGrpdCStar}, so we see that \(UV\) is an open bisection.
    Associativity of the multiplication in \(G\) then yields associativity of the composition in \(\Bis(G)\).

    For any element \(\gamma\in G\) we have \(\gamma=\gamma\gamma^{-1}\gamma\) and \(\gamma^{-1}=\gamma^{-1}\gamma\gamma^{-1}\), and so it follows quickly that \(U^{-1}=\{\gamma^{-1}:\gamma\in U\}\) gives a partial inverse for \(U\) in \(\Bis G\).
    That this is the unique element of \(\Bis G\) follows from the fact that the idempotent semilattice in \(\Bis G \) is \(E \Bis G=\Oo G^{(0)}\), where the multiplication in \(E\Bis G\) becomes the intersection in \(\Oo G^{(0)}\), which is commutative.
\end{example}

The inverse semigroup \(\Bis G\) is also closed under binary intersections, and further has some (but not all) unions.
The union \(U\cup V\) of two open bisections \(U,V\in\Bis G\) is again an open bisection if the source and range maps remain injective on the union.
Equivalently, \(U\cup V\) is an open bisection if and only if \(UV^{-1}\) and \(U^{-1}V\) consist only of units.
Moreover, if \((U_\alpha)_\alpha\) is a family of open bisections, the union \(\bigcup_\alpha U_\alpha\) is an open bisection if and only if \(U_\alpha U_\beta^{-1}\) and \(U_\alpha^{-1}U_\beta\) are contained in \(G^{(0)}\) for all \(\alpha,\beta\).
Since open subsets of the unit space of \(G^{(0)}\) are exactly the idempotents of \(\Bis G\).
This motivates the following definition.

\begin{definition}[{\cite[Section~1]{LawsonLenz_PseuGrpsEtaleGrpds}}]\label{defn-comp+pseuGrp}
    Let \(S\) be an inverse semigroup.
    Two elements \(t,u\in S\) are \emph{compatible} if \(tu^{-1},t^{-1}u\in ES\).
    A family \((t_\alpha)_\alpha\) of elements \(t_\alpha\in S\) is compatible if it is pairwise compatible.

    A \emph{pseudogroup} \(S\) is an inverse semigroup for which
    \begin{enumerate}
        \item binary meets exist: the infimum \(t\wedge u\in S\) exists for all \(t,u\in S\);
        \item compatible joins exist: if \((t_\alpha)_\alpha\) is a compatible family of elements in \(S\), then the supremum \(\bigvee_\alpha t_\alpha\) exists in \(S\);
        \item multiplication distributes over joins: if \((t_\alpha)_\alpha\) is a compatible family of elements in \(S\) then \(u\bigvee_\alpha t_\alpha=\bigvee_\alpha ut_\alpha\).
    \end{enumerate}

    A \emph{pseudogroup morphism} from \(S\) to \(T\) is a semigroup homomorphism \(\varphi\colon S\to T\) which preserves compatible joins.
\end{definition}

We may sometimes also refer to pseudogroup morphisms as \emph{pseudogroup homomorphisms}, \emph{homomorphisms of pseudogroups}, or simply \emph{homomorphisms} if the context is clear.

Pseudogroups together with pseudogroup morphisms form a category \(\Pseu\).
A frame \(F\) can be viewed as an idempotent pseudogroup with multiplication given by \(ef:=e\wedge f\) and involution \(e^*:=e\), so all elements of a frame (considered as a pseudogroup) are idempotents.
Conversely, any pseudogroup consisting only of idempotents will have arbitrary joins, since all families of elements will be compatible.
Pseudogroup morphisms between frames are then exactly frame homomorphisms; they preserve all joins since all families of elements of a frame are compatible in the sense of Definition~\ref{defn-comp+pseuGrp}, and they preserve meets since the meet is implemented by multiplication in the pseudogroup.
Thus we see that \(\Frame\) embeds into \(\Pseu\) as the full subcategory spanned by pseudogroups consisting only of idempotent elements.

The category \(\Grp\) of groups also embeds into \(\Pseu\) by sending a group \(\Gamma\) to the set \(\Gamma\sqcup\{0\}\), adjoining a zero.
The zero must be adjoined since pseudogroups require the empty supremum.
The only compatible subsets of \(\Gamma\sqcup\{0\}\) are the empty set, the singletons, and sets of the form \(\{g,0\}\) for some \(g\in \Gamma\).
Any group homomorphism \(\Gamma\to\Lambda\) then extends to a pseudogroup homomorphism by mapping zero to zero.

\begin{remark}
    We do not explicitly require pseudogroup homomorphisms to preserve binary meets.
    This is intentional as preservation of binary meets is too strong a property.
    For example, if \(\varphi\colon \Gamma\to\Lambda\) is a group homomorphism and we suppose that its extension \(\tilde\varphi\colon\Gamma\sqcup\{0\}\to\Lambda\sqcup\{0\}\) to pseudogroups preserved meets, for any \(g\in\Gamma\setminus\{1\}\) we would have \(0=\tilde\varphi(0)=\tilde\varphi(g\wedge 1_\Gamma)=\varphi(g)\wedge\varphi(1)\) implying \(\varphi(g)\neq \varphi(1)\).
    This shows that a meet-preserving homomorphism between groups is injective (and it is not difficult to see the converse also holds).
    Since we do not wish to restrict ourselves to injective group homomorphisms, we necessarily must allow pseudogroup homomorphisms which do not preserve meets.
\end{remark}

While pseudogroup homomorphisms need not preserve all binary meets, they do preserve compatible meets.

\begin{lemma}\label{lem-pseuHomPresCompMeet}
    Let \(\varphi\colon S\to T\) be a pseudogroup homomorphism and suppose that \(t,u\in S\) are compatible.
    Then \(\varphi(t\wedge u)=\varphi(t)\wedge\varphi(u)\).
    \begin{proof}
        For any compatible \(t,u\) in a pseudogroup we have the following:
        we note that \(tu^*u\leq t,u\) since both \(tu^*\) and \(u^*u\) are idempotents, so \(tu^*u\leq t\wedge u\).
        But since \((t\wedge u)^*(t\wedge u)\leq t^*t,u^*u\) we have
        \[t\wedge u=(t\wedge u)t^*tu^*u=tu^*u\wedge ut^*t\leq tu^*u,\]
        so \(tu^*u=t\wedge u\).
        Now we have \(\varphi(t\wedge u)=\varphi(tu^*u)=\varphi(t)\varphi(u)^*\varphi(u)=\varphi(t)\wedge\varphi(u)\), where the final equality holds since \(\varphi(t)\) and \(\varphi(u)\) are compatible in \(T\).
    \end{proof}
\end{lemma}

We have already observed that a pseudogroup morphism between frames is a frame homomorphism since the meet is implemented by multiplication.
Lemma~\ref{lem-pseuHomPresCompMeet} shows that this property is not exclusive to the idempotent frame of a pseudogroup, but holds on any compatible subset.

We briefly state and prove some of the fundamental properties of pseudogroups.

\begin{lemma}\label{lem-pseuProps}
    Each pseudogroup \(S\) has the following properties:
    \begin{enumerate}
        \item\label{lem-pseuProps-multOrder}  multiplication respects the order structure: if \(t\leq t'\) and \(u\leq u'\) then \(tu\leq t'u'\);
        \item\label{lem-pseuProps-joinMeetPres} joins and meets preserve the involution: if \((t_\alpha)_\alpha\) is a compatible family of elements in \(S\) then \(\left(\bigvee_\alpha t_\alpha\right)^*=\bigvee_\alpha t_\alpha^*\), and \((t\wedge u)^*=t^*\wedge u^*\) for all \(t,u\in S\);
        \item multiplication distributes over joins on the right: if \((t_\alpha)_\alpha\) is a compatible family of elements in \(S\) then \(\left(\bigvee_\alpha t_\alpha\right)u=\bigvee_\alpha (t_\alpha u)\);
        \item multiplication distributes over meets: for \(t,u,v\in S\) we have \(t(u\wedge v)=tu\wedge tv\) and \((u\wedge v)t=ut\wedge vt\).
    \end{enumerate}
    \begin{proof}
        \begin{enumerate}
            \item This holds more generally for inverse semigroups.
            By assumption we have \(t=tt^*t'\) and \(u=u'u^*u\).
            Thus we have \(tu=tt^*t'u'u^*u\), hence \(tu\leq t'u'\) by \cite[Lemma~1.4.6, Proposition~1.4.7]{Lawson_InvSemiGrps}.
            \item Since the \({}^*\)-operation is an order isomorphism by Lemma~\ref{lem-starOrdPres}, we see immediately that it preserves any suprema and infima that exist.
            \item By part~(\ref{lem-pseuProps-joinMeetPres}) we have
            \[\bigvee_\alpha(t_\alpha u)=\left(\bigvee_\alpha(t_\alpha u)^*\right)^*=\left(\bigvee_\alpha u^*t_\alpha^*\right)^*=\left(u^*\bigvee_\alpha t_\alpha^*\right)^*=\left(\bigvee_\alpha t_\alpha\right)u.\]
            \item Part~(\ref{lem-pseuProps-multOrder}) implies \(t(u\wedge v)\leq tu\wedge tv\).
            If \(w\in S\) satisfies \(w\leq tu\) and \(w\leq tv\) we have \(t^*w\leq t^*tu,t^*tv\) by part~(\ref{lem-pseuProps-multOrder}), and hence \(t^*w\leq u,v\).
            Thus \(t^*w\leq u\wedge v\), whereby \(w=tt^*w\leq t(u\wedge v)\).
            Thus \(t(u\wedge v)\) is the infimum of \(tu\) and \(tv\). \qedhere
        \end{enumerate}
    \end{proof}
\end{lemma}

\begin{lemma}\label{lem-compSameSourceEqual}
    Let \(S\) be a pseudogroup and suppose that \(t,u\in S\) are compatible.
    If \(t^*t=u^*u\) then \(t=u\).
    \begin{proof}
        The join \(t\vee u\) exists in \(S\) by assumption and is greater than each \(t\) and \(u\), so we have
        \begin{equation*}
            t=(t\vee u)t^*t=(t\vee u)u^*u=u.\qedhere
        \end{equation*}
    \end{proof}
\end{lemma}

The bisection inverse semigroup \(\Bis G\) of an \'etale groupoid \(G\) in Example~\ref{ex-bisGInvSemGrp} is also our main example of a pseudogroup.
By design, compatible families of bisections are exactly those whose unions are bisections, and one readily checks that multiplication distrubtes over these unions.

Pseudogroups are necessarily unital, the set of idempotents \(ES\) of a pseudogroup \(S\) is a compatible set, and its join functions as a unit.
Every pseudogroups also has a zero, realised by the empty join.
Unlike frames, pseudogroups do not generally have a top element with respect to the order structure.
In the example of \(\Bis G\) for an \'etale groupoid \(G\), unions of open bisections are not generally bisections, and as a result the pseudgroup \(\Bis G\) contains a top element if and only if \(G=G^{(0)}\).

Since pseudogroup morphisms preserve compatible joins, they necessarily map zeros to zeros and units to units.

Similar to many categories where the morphisms are functions which preserve an algebraic-flavoured structure, the all bijective homomorphisms are isomorphisms.

\begin{lemma}\label{lem-bijPseuHomIsIso}
    A pseudogroup homomorphism \(\varphi\colon S\to T\) is an isomorphism if it is bijective.
    \begin{proof}
        Let \(\psi\) be the inverse function to \(\varphi\).
        For \(t,u\in S\) we have 
        \[\varphi(\psi(t)\psi(u))=\varphi(\psi(t))\varphi(\psi(t))=tu=\varphi(\psi(tu)),\]
        whereby \(\psi(t)\psi(u)=\psi(tu)\) as \(\varphi\) is injective.

        For a compatible family \((t_\alpha)_\alpha\in T\) we have
        \[\varphi\left(\psi\left(\bigvee_\alpha t_\alpha\right)\right)=\bigvee_\alpha t_\alpha=\bigvee_\alpha \varphi(\psi(t_\alpha))=\varphi\left(\bigvee_\alpha \psi(t_\alpha)\right),\]
        and injectivity of \(\varphi\) shows that \(\psi\) preserves the join of the family \((t_\alpha)_\alpha\).
    \end{proof}
\end{lemma}

\subsection{Inverse semigroup actions}

Throughout we assume that the inverse semigroups considered have both a zero and a unit.
If they do not already have these, these can be easily adjoined.
We refer to \cite[Section~4]{Exel_InvSmiGrpsCombCStarAlgs} as a general reference for inverse semigroup actions and transformation groupoids.
The following constructions can be found there.

\begin{definition}
    A \emph{partial homeomorphism} of a topological space \(X\) is a homeomorphism \(U\xrightarrow{\sim}V\) between two open subsets of \(X\).
    The collection of partial homeomorphisms forms an inverse semigroup \(\Pi(X)\) where the operation is given by partial composition.
    The inverse semigroup has a unit given by the identity map on \(X\), and a zero given by the empty function.

    An \emph{action} of an inverse semigroup \(S\) on \(X\) by homeomorphisms is a semigroup homomorphism \(\alpha\colon S\to\Pi(X)\), denoted \(t\mapsto\alpha_t\), with \(\alpha_0=\emptyset\) and \(\alpha_1=\Id_X\).
\end{definition}

Given an action \(\alpha\) of an inverse semigroup \(S\) on a space \(X\), each \(t\in S\) has an open domain \(\dom(t)\subseteq X\) of points on which it acts; the domain of the partial homeomorphism \(\alpha_t\).
We often write \(t\cdot x:=\alpha_t(x)\) for points \(x\) in the domain of \(t\), analogous to the commonly used notation for group actions.

We equip the set \(D:=\{(t,x):t\in S, x\in\dom(t)\}\) with the subspace topology inherited from \(S\times X\), where we consider \(S\) with the discrete topology.
We define an equivalence relation on \(D\) by declaring \((t,x)\sim(u,y)\) if and only if \(x=y\) and there is an idempotent \(e\in E(S)\) such that \(x=y\in\dom(e)\) and \(te=ue\).

\begin{proposition}[{\cite[Section~4]{Exel_InvSmiGrpsCombCStarAlgs}}]
    The partial multiplication on \(D/\sim\) given by
    \[[t,x][u,y]=[tu,y],\]
    defined if and only if \(x=u\cdot y\) turns \(D/\sim\) into an \'etale groupoid.
    We call this groupoid the \emph{transformation groupoid} \(S\ltimes_\alpha X\) associated to the action \(\alpha\).
    The map \(x\mapsto[1_S,x]\) is a homeomorphism \(X\to(S\ltimes_\alpha X)^{(0)}\), and we identify the two under this homeomorphism.
    The source and range maps are given by 
    \[s[t,x]=x,\qquad r[t,x]=t\cdot x,\]
    and the inverse is given by
    \[[t,x]^{-1}=[t^*,t\cdot x].\]
\end{proposition}

We will often simply write \(S\ltimes X\) if the action \(\alpha\) is understood. 
Each element \(t\in S\) of the inverse semigroup specifies a subset \(U_t:=\{[t,x]:x\in\dom(t)\}\subseteq S\ltimes X\).
This is an open bisection of \(S\ltimes X\) by \cite[4.18~Proposition]{Exel_InvSmiGrpsCombCStarAlgs}.
These open sets cover \(S\ltimes X\), and the source and range maps restrict to local homeomorphisms on them.

%% file: sections/cats.tex
In this section we discuss the functors between \(\EGA\) and \(\Pseu\) which will form the adjunction in Section~\ref{sec-adjunction}.
On objects, these functors have been discussed in various places including \cite{MatsnevResende_EtaleGrpdsGermGrpds}, \cite{LawsonLenz_PseuGrpsEtaleGrpds}, and \cite{CockettGarner_GeneralisingEtaleGrpds}.
Each of these references treats a different class of morphisms, and our approach is closest to Cockett and Garner's work \cite{CockettGarner_GeneralisingEtaleGrpds}.
The construction of these functors and the following adjunction in Section~\ref{sec-adjunction} are a special case of \cite[Section~6]{CockettGarner_GeneralisingEtaleGrpds}.
Our construction of \'etale groupoids from pseudogroups is equivalent to the above mentioned references, but the particular method we employ of using transformation groupoids of pseudogroup actions is inspired by \cite{BussExelMeyer_InvSemiGrpActionsGrpdActions}.
We emphasise that the groupoid construction in \cite{BussExelMeyer_InvSemiGrpActionsGrpdActions} is not equivalent to the other constructions mentioned above for reasons we shall discuss later.

\subsection{The bisection functor}

We have seen that for any \'etale groupoid \(G\), the bisections of \(G\) form a pseudogroup \(\Bis G\).
We shall show that this assignment describes the action on objects of a functor \(\Bis\colon\EGA\to\Pseu\), and that this functor extends \(\Oo\colon\Top^\op\to\Frame\) when vieweing \(\Top^\op\) and \(\Frame\) embedded into \(\EGA\) and \(\Pseu\) respectively.

Crucial to this is \cite[Lemma~2.5]{Tay_FunctGrpdCstarAlg}, which we recall here.
\begin{lemma}[{\cite[Lemma~2.5]{Tay_FunctGrpdCstarAlg}}]\label{lem-prodBisIsBis}
    Let \(h\colon G\curvearrowright H\) be an actor between \'etale groupoids \(G\) and \(H\).
    Let \(U\in\Bis(G)\) and \(V\in\Bis(H)\) be open bisections.
    The set 
    \[U\cdot V=\{\gamma\cdot x:\gamma\in G,x\in H, s(\gamma)=\rho(x)\}\]
    is an open bisection of \(H\).
\end{lemma}

Lemma~\ref{lem-prodBisIsBis} with a fixed open bisection \(V\in\Bis H\) then yields a function \(\Bis G\to\Bis H\). 
The unit of \(\Bis H\) is the open bisection \(H^{(0)}\), and this yields a pseudogroup homomorphism.

\begin{proposition}\label{prop-FunctorBis}
    The map \(\Bis h\colon\Bis G\to\Bis H\) given by 
    \[\Bis h(U)=U\cdot_h H^{(0)}\]
    is a pseudogroup morphism. 
    The assignment \(h\mapsto\Bis h\) defines a functor 
    \[\Bis\colon \EGA\to\Pseu\]
    extending the functor \(\Oo\colon\Top^\op\to\Frame\).
    \begin{proof}
        Lemma~\ref{lem-prodBisIsBis} shows that \(\Bis h(U)\) belongs to \(\Bis H\) for any \(U\in\Bis G\), so this function is well-defined.
        For \(U,V\in\Bis G\), a generic element of \(\Bis h(U)\Bis h(V)\) has the form \((\gamma\cdot x)(\eta\cdot y)\) for some \(\gamma\in U\), \(\eta\in V\), and \(x,y\in H^{(0)}\).
        Using the axioms of Definition~\ref{defn-actor} and that \(x\) is a unit in \(H\) we see
        \begin{equation*}
            (\gamma\cdot x)(\eta\cdot y)=\gamma\cdot(x(\eta\cdot y))=\gamma\cdot(\eta\cdot y)=(\gamma\eta)\cdot y,
        \end{equation*}
        whereby \(\Bis h(U)\Bis h(V)\subseteq\Bis h(UV)\).
        A symmetric argument gives the reverse inclusion, picking \(x=r(\eta\cdot y)\) in the above setup.
        Thus \(\Bis h(U)\Bis h(V)=\Bis h(UV)\) and so \(\Bis h\) is a semigroup homomorphism.

        Let \((U_\alpha)_\alpha\subseteq\Bis G\) be a compatible family.
        Since multiplication distributes over union in \(\Bis H\), we have
        \begin{align*}
            \Bis h\left(\bigcup_\alpha U_\alpha\right)&=\left\{\gamma\cdot x: \gamma\in \bigcup_\alpha U_\alpha, x\in H^{(0)}, s(\gamma)=\rho(x)\right\}\\
            &=\bigcup_\alpha\{\gamma\cdot x: \gamma\in U_\alpha, x\in H^{(0)}, s(\gamma)=\rho(x)\}\\
            &=\bigcup_\alpha U_\alpha\cdot H^{(0)}\\
            &=\bigcup_\alpha \Bis h(U_\alpha),
        \end{align*}
        hence \(\Bis h\) preserves compatible joins.

        For two bisections \(U,V\in\Bis G\) we have
        \begin{align*}
            \Bis h(U\cap V)&=(U\cap V)\cdot H^{(0)}\\
            &=\{\gamma\cdot t:\gamma\in U\cap V, t\in H^{(0)}, s(\gamma)=\rho(t)\}\\
            &=U\cdot H^{(0)}\cap V\cdot H^{(0)}\\
            &=\Bis h(U)\cap \Bis h(V),
        \end{align*}
        so \(\Bis h\) is a pseudogroup homomorphism.

        For actors \(h\colon G\curvearrowright H\), \(k\colon H\curvearrowright K\) and an open bisection \(U\in\Bis G\) we have
        \begin{align*}
            \Bis(kh)(U)&=U\cdot_{kh}K^{(0)}\\
            &=\{\gamma\cdot_{kh}t:\gamma\in U, t\in K^{(0)}, s(\gamma)=\rho_{kh}(t)\}\\
            &=\{(\gamma\cdot_h \rho_k(t))\cdot_k t:\gamma\in U, t\in K^{(0)}, s(\gamma)=\rho_{kh}(t)\}\\
            &=(U\cdot_h H^{(0)})\cdot_k K^{(0)}\\
            &=\Bis k(\Bis h(U)),
        \end{align*}
        so \(\Bis\) preserves the composition.
        For the particular actor \(1_G\colon G\curvearrowright G\) given by left multiplication, we have \(U\cdot_{1_G}G^{(0)}=UG^{(0)}=U\), hence \(\Bis 1_G=1_{\Bis G}\), and so \(\Bis\) is a functor.

        Finally, we show that \(\Bis\) restricts to \(\Oo\) on the full subcategory \(\Top^\op\).
        Let \(X\) and \(Y\) be topological space and fix an actor \(X\curvearrowright Y\) i.e. a continuous anchor map \(\rho\colon Y\to X\). 
        For \(U\in\Bis X=\Oo X\), the frame homomorphism \(\Bis \rho\) send \(U\) to
        \[U\cdot Y=\{x\cdot y: x\in U, y\in Y, x=s(x)=\rho(y)\}.\]
        But since \(Y\) constists only of units, the definition of an actor yields \(y=s(y)=s(x\cdot y)=x\cdot y\).
        Thus \(U\cdot Y=\{y: x\in U, y\in Y, x=\rho(y)\}=\rho^{-1}(U)\), and so \(\Bis \rho\) is the preimaging function associated to \(\rho\), which is exactly \(\Oo\rho\).
    \end{proof}
\end{proposition}

Just as spatial frames are those in the essential image of \(\Oo\), we say a pseudogroup \(S\) is \emph{spatial} if \(S\) is isomorphic to \(\Bis G\) for some \'etale groupoid \(G\).
The full subcategory spanned by spatial pseudogroups is denoted by \(\Pseu_\Spat\).

\subsection{The spatialisation functor}

There are a number of existing ways which (re)construct an \'etale groupoid from pseudogroup-like data, including but not limited to \cite{Resende_EtaleGrpdsQuantales}, \cite{LawsonLenz_PseuGrpsEtaleGrpds}, \cite{BussExelMeyer_InvSemiGrpActionsGrpdActions}, \cite{CockettGarner_GeneralisingEtaleGrpds}.
We shall employ a similar but distinct construction to \cite{BussExelMeyer_InvSemiGrpActionsGrpdActions} by constructing transformation groupoids associated to an action on a character space.
The resulting groupoid we construct will be equivalent to  the constructions in \cite{Resende_EtaleGrpdsQuantales}, \cite{LawsonLenz_PseuGrpsEtaleGrpds}, and \cite{CockettGarner_GeneralisingEtaleGrpds}.

A pseudogroup \(S\) acts on the spectrum of its frame of idempotents in a canonical way.
Given an element \(t\in S\) and a character \(\chi\) on the frame \(ES\) of idempotents with \(\chi(t^*t)=1\), the map
\[e\mapsto\chi(tet^*)\]
forms a character on \(ES\).
One readily verifies that this defines a homeomorphism \(c_t\colon \Uu_{t^*t}\to\Uu_{tt^*}\), and that \(c_{tu}=c_t\circ c_u\) for all \(t,u\in S\).
Thus the homeomorphisms \((c_t)_{t\in S}\) describe an action of \(S\) on \(\widehat{ES}\).

\begin{definition}\label{defn-specPseu}
    The \emph{spectrum} of a pseudogroup \(S\) is the \'etale groupoid \(\Sigma S:=S\ltimes\widehat{ES}\)
\end{definition}

Buss, Exel, and Meyer \cite{BussExelMeyer_InvSemiGrpActionsGrpdActions} consider the action of an inverse semigroup on the space of characters of the idempotent semilattice and an inverse semigroup.
This differs from the characters we consider in Definition~\ref{defn-frameChar} as they are not required to preserve arbitrary joins (indeed, the idempotent semilattice of an inverse semigroup need to have arbitrary joins of which to speak).
The resulting groupoid constructed in \cite[Section~3]{BussExelMeyer_InvSemiGrpActionsGrpdActions} is not \(T_1\) unless it is empty or the character space contains only one element (the case where the inverse semigroup is a group). 
In particular, it rarely recovers an \'etale groupoid from its bisection inverse semigroup, which is a desirable property we will require for later applications (at least, this property will hold on an interesting subcategory of \'etale groupoids). 

The spectrum of a pseudogroup as in Definiton~\ref{defn-specPseu} is equivalent to both the constructions in \cite{MatsnevResende_EtaleGrpdsGermGrpds} and \cite{LawsonLenz_PseuGrpsEtaleGrpds}.
The difference of our approach to these is in the treatment of morphisms: Matsnev and Resende provide a duality for isomorphisms of pseudogroups, whereas Lawson and Lenz consider \emph{callitic} morphisms, which are pseudogroup morphisms with some further conditions allowing the construction of a continuous functor between the associated groupoids.
While pseudogroup morphisms do not generally induce continuous functors, they do admit a functorial construction of actors, as we now show.

A pseudogroup homomorphism \(\varphi\colon S\to T\) restricts to a frame homomorphism \(ES\to ET\) so we can use the functor \(\sigma\) from Section~\ref{sec-prelims} to acquire a continuous map \(\rho_\varphi:=\sigma(\varphi|_{ES})=\varphi|_{ES}^*\colon\widehat{ET}\to\widehat{ES}\), and we may extend this to a map \(\Sigma T\to \widehat{ES}\) by precomposing with the range map.
We then define a function \(\Sigma S\baltimes{s}{\rho_\varphi}\Sigma T\to\Sigma T\) by
\[[t,\tau]\cdot[u,\chi]=[\varphi(t)u,\chi],\]
which we claim forms an actor.

\begin{lemma}\label{lem-actorFromPseuHom}
    Let \(\varphi\colon S\to T\) be a pseudogroup morphism and set \(\rho_\varphi:=\varphi^*\circ r\colon \Sigma T\to\widehat{ES}\).
    Then \(\rho_\varphi\) is the anchor map of an actor \(\Sigma\varphi\colon\Sigma S\curvearrowright\Sigma T\) with multiplication map \(\Sigma S\baltimes{s}{\rho_\varphi}\Sigma T\to\Sigma T\) given by
    \[[t,\tau]\cdot[u,\chi]=[\varphi(t)u,\chi].\]
    \begin{proof}
        Note that if \([t,\tau]\) and \([u,\chi]\) are composable for the proposed actor multiplication we necessarily have \(\tau=s[t,\tau]=\rho_\varphi[u,\chi]=\varphi^*(c_u(\chi))\).
        We demonstrate the conditions of Definition~\ref{defn-actor} are satisfied:
        \begin{enumerate}
            \item For \([u,\chi]\in\Sigma T\) we have
            \[\rho_\varphi[u,\chi]\cdot[u,\chi]=[1_S,\varphi^*c_u(\chi)]\cdot[u,\chi]=[\varphi(1_S)u,\chi]=[u,\chi],\]
            using that \(\varphi(1_S)=1_T\) since it is a pseudogroup morphism.
            \item Fix composable \(([t,\varphi^*(c_u(\chi))],[u,\chi])\in\Sigma S\baltimes{s}{\rho_\varphi}\Sigma T\).
            For \(e\in ES\) we have 
            \begin{align*}
                r[t,\varphi^*c_u(\chi)](e)&=[c_t(\varphi^*c_u(\chi))](e)\\
                &=[c_u(\chi)]\varphi(te t^*)\\
                &=\chi(u\varphi(t)\varphi(e)(u\varphi(t))^*)\\
                &=[c_{\varphi(t)u}(\chi)](\varphi(e))\\
                &=[\varphi^*c_{\varphi(t)u}](e)\\
                &=\rho_\varphi([\varphi(t)u,\chi])(e),
            \end{align*}
            and hence \(r[t,\varphi^*c_u(\chi)]=\rho_\varphi([t,\varphi^*c_u(\chi)]\cdot[u,\chi])\).
            \item For composable \([t',\tau'],[t,\tau]\in \Sigma S\) and \([u,\chi]\in\Sigma T\) we have
            \begin{align*}
                [t',\tau']\cdot([t,\tau]\cdot [u,\chi])&=[t',\tau']\cdot[\varphi(t)u,\chi]\\
                &=[\varphi(t')\varphi(t)u,\chi]\\
                &=[\varphi(t't)u,\chi]\\
                &=[t't,\tau]\cdot[u,\chi]\\
                &=([t',\tau'][t,\tau])\cdot[u,\chi],
            \end{align*}
            as required.
            \item For composable \([t,\tau]\in\Sigma S\) and \([u,\chi]\in\Sigma T\) we have 
            \[s([t,\tau]\cdot[u,\chi])=s[\varphi(t)u,\chi]=\chi=s[u,\chi],\]
            as required.
            \item For composable \([t,\tau]\in\Sigma S\) and \([u,\chi],[u',\chi']\in\Sigma T\) we have
            \begin{align*}
                ([t,\tau]\cdot[u,\chi])[u',\chi']&=[\varphi(t)u,\chi][u',\chi']\\
                &=[\varphi(t)uu',\chi']\\
                &=[t,\tau]\cdot[uu',\chi']\\
                &=[t,\tau]\cdot([u,\chi][u',\chi']).
            \end{align*}
        \end{enumerate}
        Lastly, we show that the multiplication map is continuous. 
        Fix a net \(([t_\lambda,\tau_\lambda],[u_\lambda,\chi_\lambda])\) in \(\Sigma S\baltimes{s}{r}\Sigma T\) converging to \(([t,\tau],[u,\chi])\).
        Then \(t_\lambda\to t\) in \(S\) and \(\tau_\lambda\to\tau\) in \(\widehat{ES}\), so there is some idempotent \(e\in ES\) with \(\tau(e)=1\) and \([t_\lambda,\tau_\lambda]=[t,\tau_\lambda]\) for sufficiently large \(\lambda\).
        Similarly, \([u_\lambda,\chi_\lambda]=[u,\chi_\lambda]\) for sufficiently large \(\lambda\).
        For sufficiently large \(\lambda\) we then have
        \[[t_\lambda,\tau_\lambda]\cdot[u_\lambda,\chi_\lambda]=[t,\tau_\lambda]\cdot[u,\chi_\lambda]=[\varphi(t)u,\chi_\lambda]\]
        which converges to \([\varphi(t)u,\chi]=[t,\tau]\cdot[u,\chi]\), as required.
        Hence there is an actor \(\Sigma\varphi\) with this anchor and multiplication.
    \end{proof}
\end{lemma}

This construction of an actor from a pseudogroup morphism is functorial.

\begin{proposition}
    The assignments \(S\mapsto\Sigma S\) and \((\varphi\colon S\to T)\mapsto\Sigma\varphi\) form a functor \(\Sigma\colon\Pseu\to\EGA\).
    \begin{proof}
        Let \(1_S\) be the identity pseudgroup morphism on \(S\).
        The induced map \(\widehat{ES}\to\widehat{ES}\) is then also the identity, and so the anchor map for \(\Sigma 1_S\) is just the range map \(\Sigma S\to\widehat{ES}\).
        For any composable \(([t,\tau],[u,\chi])\in \Sigma S\baltimes{s}{r}\Sigma S\) we have
        \[[t,\tau]\cdot_{\Sigma1_S}[u,\chi]=[1_S(\tau)u,\chi]=[tu,\chi]=[t,\tau][u,\chi],\]
        Hence \(\Sigma1_S\) is the left-multiplication on \(\Sigma S\) by itself, which is the identity at \(\Sigma S\) in \(\EGA\).
        
        Let \(\varphi\colon S\to T\) and \(\psi\colon T\to R\) be pseudogroup homomorphisms.
        Restricting to idempotent frames, we have \((\psi\circ\varphi)|_{ES}=\psi|_{ET}\circ\varphi|_{ES}\) which descends to \(\sigma(\psi\circ\varphi)=\sigma(\varphi)\circ\sigma(\psi)\colon\widehat{ER}\to\widehat{ES}\), as \(\sigma\) is a functor.
        The corresponding anchor map of \(\Sigma(\psi\circ\varphi)\) is \(\rho_{\psi\circ\varphi}=(\psi\circ\varphi)|_{ES}^*\circ r=\psi|_{ET}^*\circ r\circ\varphi|_{ES}^*\circ r=\rho_\varphi\circ\rho_\psi\).
        Thus the anchor maps of \(\Sigma (\psi\circ\varphi)\) and \((\Sigma \psi)(\Sigma\varphi)\) coincide.

        For \([t,\tau]\in \Sigma S\) and \([v,\mu]\in\Sigma R\) with \(s[t,\tau]=\rho_{\psi\circ\varphi}[v,\mu]\) we have
        \begin{align*}
            [t,\tau]\cdot_{(\Sigma\psi)(\Sigma\varphi)}[v,\mu]&=([t,\tau]\cdot_{\Sigma\varphi} [1_T,\psi^*\mu])\cdot_{\Sigma\psi}[v,\mu]\\
            &=[\varphi(t),\psi^*\mu]\cdot_{\Sigma\psi}[v,\mu]\\
            &=[\psi(\varphi(t))v,\mu]\\
            &=[t,\tau]\cdot_{\Sigma(\psi\circ\varphi)}[v,\mu],
        \end{align*}
        thus \((\Sigma\psi)(\Sigma\varphi)=\Sigma(\psi\circ\varphi)\), showing that \(\Sigma\) is a functor.
    \end{proof}
\end{proposition}

%% file: sections/adjunction.tex
In this section we establish the adjunction between the categories \(\Pseu\) and \(\EGA\) via the functors \(\Sigma\) and \(\Bis\), as well as the equivalence between the full subcategories of spatial pseudogroups and sober \'etale groupoids.
This section should be compared with Cockett and Garner's article \cite{CockettGarner_GeneralisingEtaleGrpds}, specifically Section 6.
In fact, Cockett and Garner prove a family of adjunctions in a wider framework of which our adjunction is a special case.
Our discussion and proof sacrifices generality in exchange for providing a more direct and accessible argument for readers not as familiar with the language or techniques used by Cockett and Garner.

We begin by describing the natural transformations which will form the unit and counit of the adjunction between the functors \(\Sigma\) and \(\Bis\).

\begin{proposition}\label{prop-nat1toBisSigma}
    Let \(S\) be a pseudogroup.
    For each \(t\in S\), the set \(\eta_S(t):=\{[t,\chi]:\chi\in\Uu_{t^*t}\}\) is an open bisection of \(\Sigma S\).
    Moreover, the map \(\eta_S\colon S\to\Bis\Sigma S\) sending \(t\) to \(\eta_S(t)\) is a pseudogroup homomorphism, and these maps assemble into a natural transformation \(\eta\colon 1_{\Pseu}\Rightarrow \Bis\Sigma\).
    \begin{proof}
        The set \(\eta_S(t)\) is an open bisection by \cite[4.18~Proposition]{Exel_InvSmiGrpsCombCStarAlgs}.
        For \(t,u\in S\), the product of a composable pair of elements \([t,x]\in\eta_S(t)\) and \([u,y]\in\eta_S(u)\) is \([tu,y]\), which clearly belongs to \(\eta_S(tu)\).
        Conversely, any element \([tu,y]\in\eta_S(tu)\) breaks down to the product \([t,u\cdot y][u,y]\in\eta_S(t)\eta_S(u)\), so we see that \(\eta_S\) is a semigroup homomorphism.

        Let \((t_\alpha)_\alpha\subseteq S\) be a compatible family and write \(t:=\bigvee_\alpha t_\alpha\).
        Lemma~\ref{lem-invSGrpHomOrdPres} implies that \(\eta_S(t)\) is an upper bound for all the \(\eta_S(t_\alpha)\), so we have \(\bigcup_\alpha\eta_S(t_\alpha)\subseteq\eta_S(t)\).
        Fix an element \([t,\chi]\in\eta_S(t)\).
        Using Lemma~\ref{lem-starOrdPres} and part~(\ref{lem-pseuProps-multOrder}) of Lemma~\ref{lem-pseuProps}, we see that \(t^*t=\bigvee_\alpha t_\alpha^*t_\alpha\), so \(\chi\left(\bigvee_\alpha t_\alpha^*t_\alpha\right)=1\).
        Since \(\chi\) is a character, there necessarily exists \(\alpha\) such that \(\chi(t_\alpha^*t_\alpha)=1\).
        We then have \([t,\chi]=[t_\alpha,\chi]\), whereby \([t,\chi]\in \eta_S(t_\alpha)\).
        This shows the reverse inclusion, so \(\eta_S\) preserves compatible joins.
        
        Fix \(t,u\in S\).
        Again using Lemma~\ref{lem-invSGrpHomOrdPres} we see that \(\eta_S(t\wedge u)\subseteq \eta_S(t)\cap \eta_S(u)\).
        For the reverse inclusion, note that an element \(\eta_S(t)\cap\eta_S(u)\) has the form \([t,\chi]=[u,\chi]\) for some \(\chi\in\widehat{ES}\) with \(\chi(t^*t)=\chi(u^*u)=1\).
        Moreover, since the germs agree, there is an idempotent \(e\in ES\) with \(te=ue\) and \(\chi(e)=1\).
        We then gain \(te\leq t\wedge u\) \([t,\chi]=[te,\chi]=[t\wedge u,\chi]\), hence \([t,\chi]\in\eta_S(t\wedge u)\).
        Thus \(\eta_S\) preserves binary meets, and hence is a pseudogroup morphism.

        Lastly, we show naturality. 
        Fix a pseudogroup morphism \(\varphi\colon S\to T\).
        For \(t\in S\) we have
        \[\Bis\Sigma\varphi(\eta_S(t))=\eta_S\cdot_{\Sigma S}\widehat{ET}=\{[\varphi(t),\tau]:\tau\in\widehat{ES},\tau(\varphi(t^*t))=1\}=\eta_T(\varphi(t)),\]
        yielding naturality.
    \end{proof} 
\end{proposition}

If \(S=F\) is a frame (i.e. consists only of idempotents), the notation map \(\eta_F\) is now overloaded: it is defined both in Proposition~\ref{prop-nat1toBisSigma} and in Lemma~\ref{lem-frameSpecTopo}.
Since we identify \(\widehat{ES}=\sigma F\) with \((\Sigma S)^{(0)}\), and the functor \(\Bis\) restricts to \(\Oo\) on the full subcategory \(\Top^\op\) of \(\EGA\), the definition of \(\eta\) in Proposition~\ref{prop-nat1toBisSigma} extends the previously existing definition from Lemma~\ref{lem-frameSpecTopo}, so this overloading of notation is justified.

\begin{proposition}\label{prop-counitActorDefn}
    Let \(G\) be an \'etale groupoid.
    There is an actor \(\epsilon_G\colon \Sigma\Bis G\curvearrowright G\) with anchor map \(\rho_G:=\kappa_{G^{(0)}}\circ r\colon G\to(\Sigma\Bis G)^{(0)}\), where \(\kappa_{G^{(0)}}\colon G^{(0)}\to \widehat{\Oo G^{(0)}}\) is the continuous map sending a point to its detection character.
    The actor multiplication \(\cdot_{\epsilon_G}\colon\Sigma\Bis G\baltimes{s}{\rho_G} G\to G\) is given by
    \[[U,\chi]\cdot_{\epsilon_G}\gamma=U\gamma.\]
    \begin{proof}
        For brevity we write \(\kappa\) for \(\kappa_{G^{(0)}}\).
        We first note that the actor multiplication described is well-defined, since if \([U,\chi]=[V,\chi]\) for two open bisections and \(\chi=\rho_G(\gamma)=\kappa(r(\gamma))\), then there is an open subset \(W\subseteq G^{(0)}\) with \(UW=VW\) and \(\kappa(r(\gamma))(W)=1\), whereby \(r(\gamma)\in W\).
        We then have \(U\gamma=UW\gamma=VW\gamma=V\gamma\), as needed.

        Since \(\rho_G(\gamma)=[G^{(0)},\kappa_{G^{(0)}}(r(\gamma))]\) we have \(\rho_G(\gamma)\cdot\gamma=G^{(0)}\gamma=\gamma\), so the first axiom of Definition~\ref{defn-actor} is satisfied.
        We note that for \(U\in\Bis G\) with \(r(\gamma)\in s(U)\), the action of \(\Bis G\) on \(\widehat{\Oo G^{(0)}}\) sends the character \(\kappa(r(\gamma))\) to \(c_U(\kappa(r(\gamma)))=\kappa(r(\gamma))(U\cdot U^*)\), which sends \(W\in \Oo G^{(0)}\) to \(1\) exactly when \(r(\gamma)\in UWU^*\), exactly when \(r(U\gamma)=U^*r(\gamma)U\in W\).
        Thus \(r[U,\kappa(r(\gamma))]=\kappa(r(U\gamma))=\rho_G([U,\kappa(r(\gamma))]\cdot \gamma)\), yielding the second axiom.
        We also clearly have \([V,\tau]\cdot([U,\chi]\cdot\gamma)=VU\gamma=[VU,\chi]\cdot\gamma\) for all \([V,\tau],[U,\chi]\in\Sigma\Bis G\) and \(\gamma\in G\) for which the actor multiplication is defined, hence the proposed actor multiplication defines a left action of \(\Sigma\Bis G\) on \(G\).

        The final two axioms of Definition~\ref{defn-actor} follow from the fact that \(\rho_G\) factors through the range map in \(G\) by construction and the associativity of \(G\).
    \end{proof}
\end{proposition}

\begin{proposition}
    The actors \(\epsilon_G\colon \Sigma\Bis G\curvearrowright G\) of Proposition~\textup{\ref{prop-counitActorDefn}} assemble into a natural transformation \(\epsilon\colon \Sigma\Bis\Rightarrow 1_{\EGA}\).
    \begin{proof}
        Fix an actor \(h\colon G\curvearrowright H\).
        Showing naturality of \(\epsilon\) amounts to showing the equality \(\epsilon_H(\Sigma\Bis h)=h\epsilon_G\).
        The anchor map of \(\Sigma\Bis h\) is \(\rho_{\Sigma\Bis h}=(\Bis h)|_{E\Bis G}^*\circ r\), where \((\Bis h)|_{E\Bis G}^*\) sends a character \(\chi\in\widehat{\Oo H^{(0)}}\) to \(\chi\circ(\Bis h)|_{E\Bis G}\).
        For \(x\in H\), the anchor of \(\epsilon_H(\Sigma\Bis h)\) sends \(x\) to the character
        \[U\mapsto [\kappa_{H^{(0)}}(r(x))](\Bis h(U))=[\kappa_{H^{(0)}}(r(x))](\rho_h^{-1}(U))=\begin{cases}
            1,&r(x)\in\rho_h^{-1}(U),\\
            0,&r(x)\notin\rho_h^{-1}(U).
        \end{cases}\]
        The anchor of \(h\epsilon_G\) sends \(x\) to \(\kappa_{G^{(0)}}(\rho_h(r(x)))\).
        This character is defined by
        \[[\kappa_{G^{(0)}}(\rho_h(r(x)))](U)=\begin{cases}
            1,&\rho_h(r(x))\in U,\\
            0,&\rho_h(r(x))\notin U,
        \end{cases}=\begin{cases}
            1,&r(x)\in\rho_h^{-1}(U),\\
            0,&r(x)\notin\rho_h^{-1}(U),
        \end{cases}\]
        which is exactly \([\kappa_{H^{(0)}}(r(x))](\Bis h(U))\).
        Hence \(\rho_{\epsilon_H(\Sigma\Bis h)}=\rho_{h\epsilon_G}\).

        The pair \(([U,\chi],x)\in \Sigma\Bis G\times H\) is composable with respect to the actors of interest if and only if \(\chi=\kappa_H^{(0)}(r(x))\).
        The actor multiplication of \(h\epsilon_G\) maps this pair to 
        \[[U,\chi]\cdot_{\epsilon_G h}x=([U,\chi]\cdot_{\epsilon_G}\rho_h(x))\cdot_h x=(U\rho_h(x))\cdot_h x,\]
        and the actor multiplication of \(\epsilon_H(\Sigma\Bis h)\) maps it to
        \begin{align*}
            [U,\chi]\cdot_{\epsilon_H(\Sigma\Bis h)}x&=([U,\chi]\cdot_{\Sigma\Bis h}\kappa_{H^{(0)}}(r(x)))\cdot_{\epsilon_H} x\\
            &=[\Bis h(U),\kappa_{H^{(0)}}(r(x))]\cdot_{\epsilon_H}x\\
            &=\Bis h(U)x.
        \end{align*}
        But \(\Bis h(U)x=(\Bis h(U)r(x))x\) is the product of the unique element of \(\Bis h(U)\) with source \(r(x)\) with \(x\).
        Since \(\Bis h(U)=U\cdot_h H^{(0)}\), we have that \(\Bis h(U)r(x)=\gamma\cdot_h r(x)\), where \(\gamma\in U\) is the unique element with \(s(\gamma)=\rho_h(x)\).
        Thus \(\Bis h(U)r(x)=(U\rho_h(x))\cdot_h r(x)\) and hence
        \[[U,\chi]\cdot_{\epsilon_G h}x=(U\rho_h(x))\cdot_h x=\Bis h(U)x=[U,\chi]\cdot_{\epsilon_H(\Sigma\Bis h)}x,\]
        as required.
    \end{proof}
\end{proposition}

We now show that \(\eta\) and \(\epsilon\) form the unit and counit of an adjunction.

\begin{theorem}[{cf. \cite[Theorem~6.8]{CockettGarner_GeneralisingEtaleGrpds}}]\label{thm-adjunction}
    The functors \(\Sigma\colon\Pseu\to\EGA\) and \(\Bis\colon\EGA\to\Pseu\) form an adjunction \(\Sigma\dashv\Bis\).
    The natural transformations \(\eta\) and \(\epsilon\) of Propositions~\textup{\ref{prop-nat1toBisSigma}} and \textup{\ref{prop-counitActorDefn}} are the respective unit and counit of this adjunction.
    \begin{proof}
        We show that the triangle identities for \(\eta\) and \(\epsilon\) are satisfied.

        For a pseudogroup \(S\), the anchor map of the actor \(((\epsilon \Sigma)(\Sigma\eta))_S=\epsilon_{\Sigma S}(\Sigma\eta_S)\) is given by the composition \(\chi\mapsto\kappa_{(\Sigma S)^{(0)}}(\chi)\mapsto[\kappa_{(\Sigma S)^{(0)}}(\chi)](\eta_S(\cdot))\).
        This character acts on \(ES\) by sending \(e\) to
        \[[\kappa_{(\Sigma S)^{(0)}}(\chi)](\eta_S(e))=\begin{cases}
            1,&\chi\in \eta_S(e),\\
            0,&\chi\notin\eta_S(e),
        \end{cases}=\begin{cases}
            1,&\chi(e)=1,\\
            0,&\chi(e)=0,
        \end{cases}=\chi(e).\]
        Hence \(\rho_{\epsilon_{\Sigma S}(\Sigma\eta_S)}(\chi)=\eta_S^*\kappa_{(\Sigma S)^{(0)}}(\chi)=\chi\), so the anchor map of \(\epsilon_{\Sigma S}(\Sigma\eta_S)\) is the identity on the unit space of \(\Sigma S\).
        Since both anchor maps factor through the range map, we see that the anchor map is simply the range map on \(\Sigma S\).
        
        For composable \([t,\chi],[u,\tau]\in\Sigma S\), the actor product of the two is
        \begin{align*}
            [t,\chi]\cdot_{\epsilon_{\Sigma S}(\Sigma\eta_S)}[u,\tau]&=([t,\chi]\cdot_{\Sigma \eta_S}[1_{\Bis\Sigma S},\kappa_{(\Sigma S)^{(0)}}(r[u,\tau])])\cdot_{\epsilon_{\Sigma S}}[u,\tau]\\
            &=[\eta_S(t),\kappa_{(\Sigma S)^{(0)}}(r[V,\tau])]\cdot_{\epsilon_{\Sigma S}}[u,\tau]\\
            &=\eta_S(t)[u,\tau]\\
            &=[tu,\tau],
        \end{align*}
        which is exactly the product of \([t,\chi]\) and \([u,\tau]\) in the groupoid \(\Sigma S\).
        Thus \(\epsilon_{\Sigma S}(\Sigma\eta_S)=1_{\Sigma S}\).

        Let \(G\) be an \'etale groupoid and consider the component 
        \[((\Bis\epsilon)(\eta \Bis))_G=\Bis\epsilon_{G}(\eta_{\Bis G})\colon \Bis G\to\Bis G.\]
        For an open bisection \(U\in\Bis G\), recall that \(\eta_{\Bis G}(U)\) is the set
        \[\{[U,\chi]:\chi\in\widehat{\Oo G^{(0)}}, \chi(U^*U)=1\}.\]
        For a point \(x\in G^{(0)}\) with \(\kappa_{G^{(0)}}(x)\in U^*U\), the actor product \([U,\kappa_{G^{(0)}}]\cdot_{\epsilon_G}x\) evaluates to \(Ux\), which is the unique element of \(U\) with source \(x\).
        We also note that \(x\in G^{(0)}\) satisfies \(\rho_{\epsilon_G}(x)=[U,\chi]\) for some \(\chi\in\widehat{\Oo G^{(0)}}\) if and only if \(\chi=\kappa_{G^{(0)}}(x)\) and \([\kappa_{G^{(0)}}(x)](U^*U)=1\), if and only if \(x\in U^*U\).
        Hence we compute 
        \[\Bis\epsilon_{G}(\eta_{\Bis G}(U))=\eta_{\Bis G}(U)\cdot_{\epsilon_G} G^{(0)}=\{[U,\chi]\cdot_{\epsilon_G}x: x\in U^*U\}=U,\]
        whereby \(\Bis\epsilon_G(\eta_{\Bis G})=1_{\Bis G}\), as required.
    \end{proof}
\end{theorem}

This adjunction restricts to an equivalence between the full subcategory of sober \'etale groupoids in \(\EGA\) and the full subcategory of spatial pseudogroups in \(\Pseu\).

\begin{lemma}\label{lem-spatSob}
    Let \(S\) be a pseudogroup.
    Then \(\Sigma S\) is sober.
    \begin{proof}
        Since \(\Sigma S\) is \'etale, it is locally homeomorphic to its unit space \(\widehat{ES}\).
        Since \(ES\) is a frame, its spectrum \(\widehat{ES}\) is sober by Proposition~\ref{prop-specOfFrmSob}.
        Proposition~\ref{prop-sobLocHomeo} then implies that \(\Sigma S\) is sober.
    \end{proof}
\end{lemma}

We also recall here that Proposition~\ref{prop-sobLocHomeo} implies an \'etale groupoid is sober if and only if its unit space is.

\begin{proposition}\label{prop-counitActSoberIso}
    Let \(G\) be a sober groupoid.
    Then the counit actor \(\epsilon_G\colon\Sigma\Bis G\curvearrowright G\) is an isomorphism.
    \begin{proof}
        The unit space \((\Sigma\Bis G)^{(0)}=\widehat{\Oo G^{(0)}}\) is sober by Proposition~\ref{prop-specOfFrmSob}, so the frame homomorphism \(\kappa_{G^{(0)}}\colon G^{(0)}\to \widehat{\Oo G^{(0)}}\) is an isomorphism by Theorem~\ref{thm-soberIffSpecHomeo}.
        Hence every character on \(\Oo G^{(0)}\) is a point detection character \(\chi_x:=\kappa_{G^{(0)}}(t)\) for some \(x\in G^{(0)}\).
        Define \(\rho\colon \Sigma\Bis G\to G^{(0)}\) by \(\rho:=\kappa_{G^{(0)}}^{-1}\circ r\), so that \(\rho[U,\chi_x]=U^*xU\).
        We then construct an actor multiplication for an actor \(h\colon G\curvearrowright \Sigma\Bis G\) with anchor \(\rho\) by defining
        \[\gamma\cdot_h [U,\chi_x]=[V_\gamma U,\chi_x]\]
        for \(\gamma\in G\) and \([U,\chi_x]\in \Sigma\Bis G\) with \(s(\gamma)=\rho[U,\chi_x]=UxU^*\), where \(V_\gamma\in\Bis G\) is a bisection neighbourhood of \(\gamma\).
        This formula is independent of choice of \(V_\gamma\), since for any two \(V_\gamma,V_\gamma'\in\Bis G\) containing \(\gamma\), the germ of \(V_\gamma U\) and \(V_\gamma'U\) at \(x\) both coincide with the germ of \((V_\gamma\cap V_\gamma')U\) at \(x\).
        One readily checks that this forms an actor \(h\).

        We claim that \(h\) is the inverse actor to \(\epsilon_{G}\).
        By construction the anchor maps are inverse to each other on the unit spaces, so their compositions yield the range maps on their respective domains.
        First we compute the composition \(\epsilon_G h\).
        For composable arrows \(\gamma,\eta\in G\), let \(V_\gamma\) be a bisection containing \(\gamma\). 
        We compute
        \begin{align*}
            \gamma\cdot_{\epsilon_G h}\eta&=\gamma\cdot_h([G^{(0)},\chi_{r(\eta)}])\cdot_{\epsilon_G}\eta\\
            &=[V_\gamma,\chi_{r(\eta)}]\cdot_{\epsilon_G}\eta\\
            &=V_\gamma\eta\\
            &=\gamma\eta,
        \end{align*}
        as required.
        All that remains is to check the reverse composition.
        For composable \([U,\chi_x],[V,\chi_y]\in\Sigma\Bis G\), where \(x,y\in G^{(0)}\) are units, we have
        \begin{align*}
            [U,\chi_x]\cdot_{h \epsilon_G}[V,\chi_y]&=([U,\chi_x]\cdot_{\epsilon_G}VyV^*)\cdot_h [V,\chi_y]\\
            &=UVyV^*\cdot_h [V,\chi_y]\\
            &=[UV,\chi_y]\\
            &=[U,\chi_x][V,\chi_y],
        \end{align*}
        as required.
    \end{proof}
\end{proposition}

Two groupoids are isomorphic in \(\EGA\) if and only if they are isomorphic as topological categories (i.e. there are continuous mutually inverse functors between them) by \cite[Proposition~4.19]{MeyerZhu_Groupoids}.
Hence Proposition~\ref{prop-counitActSoberIso} shows a sober \'etale groupoid \(G\) is isomorphic to \(\Sigma\Bis G\) in the `usual sense'.

Similarly to the case for frames, we may show that a pseudogroup \(S\) is spatial if and only if it is isomorphic to the pseudogroup of bisections of its associated groupoid \(\Sigma S\).

\begin{proposition}\label{prop-pseuSpatIffBisSig}
    A pseudogroup \(S\) is spatial if and only if if and only if \(ES\) is spatial.
    In this case, the homomorphism \(\eta_S\colon S\to\Bis\Sigma S\) is an isomorphism.
    \begin{proof}
        Suppose \(S\) is spatial.
        By definition there is an \'etale groupoid \(G\) with \(S\cong \Bis G\).
        This implies that \(ES\cong E\Bis G=\Oo G^{(0)}\) is a spatial frame, and hence \(ES\cong \Oo\widehat{ES}\) by \cite[5.1~Proposition]{PicadoPultr_FramesLocales}.
        Conversely, suppose \(ES\) is spatial.
        Let \(\eta_S\colon S\to\Bis\Sigma S\) be the component of the unit in Proposition~\ref{prop-nat1toBisSigma}.
        Fix \(t,u\in S\) with \(\eta_S(t)=\eta_S(u)\).
        Then \(\eta_S(t^*t)=\eta_S(u^*u)\) and hence \(t^*t=u^*u\) as \(ES\) is spatial and \(\eta\) extends the unit of the adjunction in Theorem~\ref{thm-soberSpatAdj}. 
        Further, the germs of \(t\) and \(u\) agree at all characters in \(t^*t\), yielding \(t=tt^*t=ut^*t=u\).
        
        For surjectivity, we employ a similar strategy as in Proposition~\ref{prop-pseuSpatIffBisSig}. 
        Fix an open bisection \(U\in\Bis\Sigma S\) and write it as a union \[U=\bigvee_\alpha \eta_S(t_\alpha)\cap U\]
        where \((t_\alpha)_\alpha\subseteq S\) is a family of elements whose images under \(\eta_S\) jointly cover \(U\).
        Since \((\eta_S(t_\alpha)\cap U)^*(\eta_S(t_\alpha)\cap U)\) is an open subset of \(\widehat{ES}\) and \(ES\) is spatial, there is an element \(e_\alpha\in ES\) with \(\eta_S(e_\alpha)=(\eta_S(t_\alpha)\cap U)^*(\eta_S(t_\alpha)\cap U)\), and we get \(\eta_S(t_\alpha)\cap U=\eta_S(t_\alpha e_\alpha)\).
        The union of the \(\eta_S(t_\alpha e_\alpha)\) is an open bisection, whereby the family \((\eta_S(t_\alpha e_\alpha)_\alpha)\) is compatible, and hence so is \((t_\alpha e_\alpha)_\alpha\) since \(\eta_S\) is injective. 
        We then have
        \[U=\bigcup_\alpha \eta_S(t_\alpha e_\alpha)=\eta_S\left(\bigvee_\alpha t_\alpha e_\alpha\right),\]
        yielding surjectivity.
        Lemma~\ref{lem-bijPseuHomIsIso} finishes the proof.

        Let \(\eta_S\colon S\to\Bis\Sigma S\) be the unit morphism from Proposition~\ref{prop-nat1toBisSigma}.
        For \(t,u\in S\) with \(\eta_S(t)=\eta_S(u)\) we have \(\eta_S(t^*t)=\eta_S(u^*u)\), which implies \(t^*t=u^*u\) as \(ES\) is spatial and \(\eta_S\) restricts to an isomorphism for a spatial frame.
        For each \(\chi\in\eta_S(t^*t)=\eta_S(u^*u)\subseteq\widehat{ES}\), let \(e_\chi\in ES\) be an idempotent with \(\chi(e_\chi)=1\) and \(te=ue\).
        These exist since \([t,\chi]=[u,\chi]\) for such \(\chi\); these are exactly the elements of \(\eta_S(t)=\eta_S(u)\).
        We then have \(t^*t=u^*u\leq\bigvee_{\chi}e_\chi\), since the open sets corresponding to the \(e_\chi\) in \(\Oo\widehat{ES}\) contain \(\chi\).
        Thus we have \(t=t\bigvee_{\chi}e_\chi=\bigvee_\chi te_\chi=\bigvee_\chi ue_\chi=u\), and \(\eta_S\) is injective.

        Fix an open bisection \(U\in\Bis\Sigma S\).
        Let \(I\subseteq S\) be the set of elements \(t\in S\) with \(\eta_S(t)\cap U\neq\emptyset\).
        Since the image of \(\eta_S\) covers \(\Sigma S\), we have that
        \[U=\bigcup_{t\in I}\eta_S(t)\cap U.\]
        Moreover, \(\eta_S(t)\cap U=\eta_S(t)(\eta_S(t)\cap U)^*(\eta_S(t)\cap U)\), as \(\eta_S(t)\) and \(\eta_S(t)\cap U\) are both open bisections.
        The set \((\eta_S(t)\cap U)^*(\eta_S(t)\cap U)\) is an open subset of \(\widehat{ES}\), so there is an idempotent \(e_t\in ES\) with \((\eta_S(t)\cap U)^*(\eta_S(t)\cap U)=\eta_S(e_t)\), since \(ES\) is a spatial frame.
        Hence \(\eta_S(t)\cap U=\eta_S(t)\eta_S(e_t)=\eta_S(te_t)\), and hence
        \[U=\bigcup_{t\in I}\eta_S(t)\cap U=\bigcup_{t\in I}\eta_S(te_t).\]
        The union of open bisections in an \'etale groupoid is again an open bisection if and only if the family is compatible, so the family \((\eta_S(te_t))_{t\in I}\) is a compatible family in \(\Bis\Sigma S\).
        Since \(\eta_S\) is injective, we see that \((te_t)_{t\in I}\) is compatible in \(S\), and hence its join exists.
        We then have
        \[\eta_S\left(\bigvee_{t\in I}te_t\right)=\bigcup_{t\in I}\eta_S(te_t)=U,\]
        whereby \(\eta_S\) is surjective.
        Hence \(\eta_S\) is an isomorphism by Lemma~\ref{lem-bijPseuHomIsIso}.
    \end{proof}
\end{proposition}

We can now combine the previous results to show that the adjunction between \(\Sigma\) and \(\Bis\) restricts to an equivalence between sober groupoids and spatial pseudogroups.

\begin{corollary}\label{cor-catEquivSobSpat}
    The adjunction \(\Sigma\dashv\Bis\) restricts to an equivalence between the full subcategory of \(\EGA\) spanned by sober groupoids and the full subcategory of \(\Pseu\) spanned by spatial pseudogroups.
    \begin{proof}
        Propositions~\ref{prop-counitActSoberIso} and \ref{prop-pseuSpatIffBisSig} show the unit and counit components are isomorphisms on these subcategories, hence the subcategories are equivalent.
    \end{proof}
\end{corollary}

\begin{corollary}\label{cor-coRefSubcats}
    The full subcategory of spatial pseudogroups is a reflective subcategory of \(\Pseu\) with reflector functor \(\Bis\Sigma\).
    The full subcategory of sober groupoids is a coreflective subcategory of \(\EGA\) with coreflector functor \(\Sigma\Bis\).
\end{corollary}

%% file: sections/lims.tex
In this section we use the adjunction in Theorem~\ref{thm-adjunction} and the equivalence between sober \'etale groupoids and spatial pseudogroups of Corollary~\ref{cor-catEquivSobSpat} to study limits in the respective categories.
Limits in \(\Pseu\) may be computed directly using the forgetful functor to \(\Set\), and we use this to deduce completeness of \(\Pseu_\Spat\) and hence \(\EGA_\Sob\).
We then describe some of the main examples of limits in \(\EGA_\Sob\).
Further we describe how limits and colimits of groups behave in \(\Pseu\) and \(\EGA\), and comment on the more difficult question of colimits.

\begin{theorem}\label{thm-pseuComplete}
    The forgetful functor \(U\colon\Pseu\to\Set\) creates limits.
    In particular, \(\Pseu\) is complete.
    \begin{proof}
        Let \(F\colon J\to\Pseu\) be a small diagram.
        The limit of \(UF\) in set can be realised as the set \(\Cone(1,UF)\), where \(1\in\Set\) is a singleton.
        The universal cone \(\lambda\colon\Cone(1,UF)\Rightarrow F\) has components \(\lambda_j\colon \Cone(1,UF)\to F_j\) sending a cone \(\mu\in\Cone(1,UF)\) over \(UF\) to \(\mu_j(\ast)\), where \(\ast\in 1\) is the unique element.

        We equip \(\Cone(1,UF)\) with a `pointwise' pseudogroup structure.
        The operations will all be pointwise: for \(j\in J^0\)
        \begin{itemize}
            \item the product of two cones \(\mu^1,\mu^2\in\Cone(1,UF)\) is given by \((\mu^1\mu^2)_j(\ast)=\mu^1_j(\ast)\mu^2_j(\ast)\),
            \item the \({}^*\)-operation is given by \((\mu^*)_j(\ast)=\mu_j(\ast)^*\),
            \item the meet of two cones \(\mu^1,\mu^2\in\Cone(1,UF)\) is given by \((\mu^1\wedge\mu^2)_j(\ast)=\mu^1_j(\ast)\wedge\mu^2_j(\ast)\)
            \item a family \((\mu^\alpha)_\alpha\subseteq\Cone(1,UF)\) is compatible if and only if \((\mu^\alpha_j(\ast))_\alpha\) is a compatible family in the pseudogroup \(Fj\), so these pointwise joins exist.
            The join of the family in \(\Cone(1,UF)\) is given by \(\left(\bigvee_\alpha \mu^\alpha\right)_j(\ast)=\bigvee_\alpha \mu^\alpha_j(\ast)\).
        \end{itemize}
        One readily verifies that this structure turns \(\Cone(1,UF)\) into a pseudogroup, and the components \(\lambda_j\) of the universal cone in \(\Set\) are pseudogroup homomorphisms by construction.
        Since the forgetful functor to \(\Set\) is faithful, these components lift to a cone over \(F\) in \(\Pseu\).
        Moreover, any other cone \(\mu\colon S\Rightarrow F\) over \(F\) descends to a cone over \(UF\) and we get a function \(\tilde\mu\colon US\to\Cone(1,UF)\) such that \(\lambda_j\tilde\mu=\mu_j\) for each \(j\in J^0\).
        Since all the pseudogroup operations in \(\Cone(1,UF)\) are pointwise and the \(\lambda_j\) are evaluation maps, we see that \(\tilde\mu\) is a pseudogroup homomorphism precisely when \(\lambda_j\tilde\mu\) is for each \(j\).
        But these are exactly \(\mu_j\), which are pseudogroup homomorphisms by assumption.
        The map \(\tilde\mu\) is also unique for this property, since another such morphism would descend to a function \(US\to\Cone(1,UF)\) in \(\Set\), whereby uniqueness of such a map there and faithfulness of \(U\) ensure that it was \(\tilde\mu\) already.
    \end{proof}
\end{theorem}

Since the subcategory \(\Pseu_\Spat\) of spatial pseudogroups is reflective in \(\Pseu\) by Corollary~\ref{cor-coRefSubcats}, we may use Theorem~\ref{thm-pseuComplete} to deduce completeness of \(\Pseu_\Spat\).
Corollary~\ref{cor-catEquivSobSpat} then transfers this to \(\EGA_\Sob\).

\begin{theorem}\label{thm-spatSobComplete}
    The forgetful functor \(\Pseu_\Spat\to\Set\) creates limits.
    In particular, \(\Pseu_\Spat\) and \(\EGA_\Sob\) are complete.
    \begin{proof}
        The subcategory \(\Pseu_\Spat\) is a reflective subcategory of \(\Pseu\), so the inclusion \(\Pseu_\Spat\hookrightarrow\Pseu\) creates limits by \cite[Proposition~4.5.15]{Riehl_CategoryTheoryContext}.
        Theorem~\ref{thm-pseuComplete} shows that \(\Pseu\) is complete, hence \(\Pseu_\Spat\) is complete, and moreover the limit is constructed in \(\Pseu\) so is created by the forgetful functor to \(\Set\).

        \(\EGA_\Sob\) is equivalent to \(\Pseu_\Spat\) by Corollary~\ref{cor-catEquivSobSpat}, so completeness passes through.
    \end{proof}
\end{theorem}

\begin{remark}
    The forgetful functor \(\Pseu\to\Set\) admits a left adjoint i.e. a `free pseudogroup' functor, so itself is a right adjoint and must therefore at least preserve limits.
    Moreover this demonstrates that \(\Pseu\) has colimits of all free objects, as the left adjoint to the forgetful functor must preserve them coming from \(\Set\).
    We refer the reader to \cite[Theorem~1.4.23, Theorem~1.4.24, Theorem~6.1.2]{Lawson_InvSemiGrps} for details.
\end{remark}

Since any pseudogroup morphism necessarily maps idempotents to idempotents, restricting to the frame of idempotents yields a functor \(E\colon\Pseu\to\Frame\).
Conversely we may view any frame as a pseudogroup consisting of idempotents, and pseudogroup homorphisms between two frames are exactly frame homomorphisms.
We denote the inclusion functor by \(\iota\colon\Frame\to\Pseu\).
This functor is right-adjoint to the inclusion of the category of frames into the category of pseudogroups, yielding a coreflective subcategory.

\begin{proposition}\label{prop-frmCorefPseu}
    \(E\) is right-adjoint to the inclusion \(\iota\colon\Frame\to\Pseu\), hence \(\Frame\) is a coreflective subcategory of \(\Pseu\).
    For a frame \(F\) and a pseudogroup \(S\), the natural bijections \(\Pseu(\iota F,S)\xrightarrow{\sim}\Frame(F,ES)\) send a pseudogroup homomorphism \(\varphi^\sharp\colon\iota F\to S\) to its corestriction \(\varphi^\flat\colon F\to ES\).
    \begin{proof}
        Since \(\iota F\) consists only of idempotents, the image of \(\varphi^\sharp\colon \iota F\to S\) must be contained in \(ES\), and so the corestriction \(\varphi^\flat\colon F\to ES\), given by \(\varphi^\flat(e)=\varphi^\sharp(e)\) is well-defined. 
        Moreover, every frame homomorphism \(F\to ES\) is of this form, so the mapping \(\varphi^\sharp\mapsto\varphi^\flat\) is a bijection.
        Naturality is elementary to show.
    \end{proof}
\end{proposition}

By Proposition~\ref{prop-pseuSpatIffBisSig}, the spatiality of a pseudogroup depends only on its idempotent frame and hence limit groupoids in \(\EGA_\Sob\) can be computed in \(\Top_\Sob\).

\begin{corollary}\label{cor-EPresLims}
    \(E\) preserves limits.
    In particular, the unit space of the limit of a diagram \(F\colon J\to\EGA_\Sob\) is given by the corresponding limit in \(\Top^\op_\Sob\) of the (opposite) anchor maps.
\end{corollary}

\subsection{Examples of limits}\label{subsec-EGALims}

With our newfound powers we can better describe some of the limits in \(\EGA_\Sob\).
The composition \(\EGA_\Sob\xrightarrow\Pseu\to\Set\) now gives a functor from the cateogry of sober \'etale groupoids to the category of sets which creates the desired limits, which means limits of sober groupoids in this category can be constructed in the category of sets and then `lifted' through this functor.
However, since this functor does not send a groupoid to the underlying set of arrows in the groupoid (but rather to the set of bisections), the resulting limit constructions may not be so familiar.
One example is the product in \(\EGA_\Sob\) as it does not coincide with the Cartesian product.
The product can be described even in the larger category \(\EGA\) without needing to directly invoke Theorems~\ref{thm-pseuComplete} or \ref{thm-spatSobComplete}.

\begin{proposition}\label{prop-productDisjUn}
    The product in the category \(\EGA\) is given by disjoint union.
    \begin{proof}
        Let \((G_\alpha)_{\alpha\in I}\) be a family of \'etale groupoids.
        Let \(G:=\bigsqcup_{\alpha\in I}G_\alpha\) be the disjoint union, equipped with the disjoint union topology and disjoint union groupoid structure (that is, the coproduct in the category of groupoids with functors as morphisms).
        For each \(\alpha\in I\), let \(\iota_\alpha\colon G_\alpha\to G\) be the inclusion map.

        We define a universal cone over the diagram \((G_\alpha)_\alpha\) as follows: for each \(\alpha\in I\) let \(\rho_\alpha\colon G_\alpha\to G^{(0)}\) be the composition \(\rho_\alpha:=\iota_\alpha\circ r_\alpha\), where \(r_\alpha\colon G_\alpha\to G_\alpha^{(0)}\) is the range map in \(G_\alpha\).
        For \(\gamma\in G\) and \(\eta\in G_\alpha\), we have \(s(\gamma)=\rho_\alpha(\eta)=\iota_\alpha(r_\alpha(\eta))\) only if \(\gamma\) belongs to the \(G_\alpha\)-component of the disjoint union.
        Identifying each \(G_\alpha\) with the image \(\iota_\alpha(G_\alpha)\subseteq G\), we see that \((\gamma,\eta)\) is a composable pair in \(G_\alpha\), and we can then define \(\cdot_\alpha\colon G\baltimes{s}{\rho_\alpha}G_\alpha\to G_\alpha\) by
        \[\gamma\cdot_\alpha\eta:=\gamma\eta,\]
        using the multiplication in \(G_\alpha\).

        One readily verifies that this defines a family of actors \(\pi_\alpha\colon G\curvearrowright G_\alpha\).
        Given another family of actors \((h_\alpha\colon H\curvearrowright G_\alpha)_\alpha\), we construct an actor \(h\colon H\curvearrowright G\) with \(h_\alpha=\pi_\alpha h\) for all \(\alpha\).
        The anchor map of \(h\) is given by taking the disjoint union of the anchors \(\rho_{h_\alpha}\colon G_\alpha\to H^{(0)}\) to a continuous map \(\rho_h:=\bigsqcup \rho_{h_\alpha}\colon G\to H\).
        For a pair \(x\in H\) and \(\gamma\in G\) satisfying \(s(x)=\rho_h(\gamma)\), there is a unique \(\alpha\) such that \(\gamma\in G_\alpha\) and \(\rho_h(\gamma)=\rho_{h_\alpha}(\gamma)\).
        For such a pair, we define \(x\cdot_h\gamma:=x\cdot_{h_\alpha}\gamma\).
        This specifies an actor, and by construction we have \(h_\alpha=\pi_\alpha h\).
        It is also clear that \(h\) is the unique such actor, since acting on each \(G_\alpha\) is specified and these subsets cover \(G\).
    \end{proof}
\end{proposition}

The construction of the equaliser is a little more opaque.

Let \(G\) and \(H\) be sober \'etale groupoids and let \(h,k\colon G\curvearrowright H\) be actors.
To compute the equaliser, we consider the induced homomorphisms \(\Bis h,\Bis k\colon\Bis G\to\Bis H\) and take the equaliser there.
Using Theorem~\ref{thm-spatSobComplete}, the equaliser of these two homomorphisms is the subpseudogroup \(S:=\{U\in\Bis G: \Bis h(U)=\Bis k(U)\}\subseteq\Bis G\).
Furthermore, \(S\) is spatial, so \(S\cong\Bis\Sigma S\) together with the induced homomomorphism \(\Bis\Sigma S\to \Bis G\) forms the equaliser in \(\Pseu_\Spat\).

Pushing this through \(\Sigma\) yields the groupoids \(S\ltimes\widehat{ES}\), so we wish to describe the action of \(S\) on \(\widehat{ES}\).
Firstly, we note that \(ES\) is the subset of \(\Oo G^{(0)}\) on which \(\Bis h\) and \(\Bis k\) agree, as \(E\) preserves limits by Proposition~\ref{prop-frmCorefPseu}.
Passing through the equivalence in Theorem~\ref{thm-soberSpatAdj} from spatial frames to the opposite category of sober topological spaces, we see that the unit space of \(\Sigma S\) is given by the coequaliser of the induced continuous maps \(H^{(0)}\to G^{(0)}\) from \(\Bis h\) and \(\Bis k\).
But these are exactly the anchor maps restricted to the unit space of \(H\), so \(\widehat{ES}\cong\coeq(\rho_h,\rho_k)\), which is the quotient of \(G^{(0)}\) by the equivalence relation generated by \(\rho_h(x)\sim\rho_k(x)\) for all \(x\in H\).

Viewing a unit of \(\Sigma S\) as an equivalence class \([t]\) for some unit \(t\in G^{(0)}\), the action of \(U\in S\) on \([t]\) is defined if \(U^*U\) contains a representative of \([t]\), and the action sends \([t]\) to \([UtU^*]\).
The spatialisation of \(S\) consists of such pairs \([U,[t]]\).
The actor \(\Sigma S\curvearrowright G\) has anchor map given by the quotient map \(G\to\coeq(\rho_h,\rho_k)\) sending \(\gamma\) to the class \([r(\gamma)]\).
The actor multiplication is then given by 
\[[U,[r(\gamma)]]\cdot \gamma=U\gamma.\]

Another way to view the coequaliser is by taking the open subgroupoid \(K:=\bigcup S\subseteq G\) covered by \(S\).
We place an equivalence relation on \(K\) extending the coequaliser relation for \(\rho_h\) and \(\rho_k\) by declaring \(\gamma\sim\eta\) if \(s(\gamma)\sim s(\eta)\) and there is an element \(U\in S\) containing both \(\gamma\) and \(\eta\).
The quotient of \(K\) by this relation then yields a groupoid isomorphic to \(\Sigma S\), and the actor multiplication is characterised for \([\gamma]\in K/\sim\) and \(\eta\in G\) with \(s[\gamma]=[r(\eta)]\) by \([\gamma]\cdot\eta=\gamma'\eta\), where \(\gamma'\) is the unique element of the class \([\gamma]\) which is composable with \(\eta\).

Computing pullbacks combines these two previous descriptions.

Let \(h_i\colon G_i\curvearrowright H\) be actors for \(i=1,2\) and let \(h_i\pi_i\) be the compositions with the universal product actors \(\pi_i\colon G_1\times G_2\curvearrowright G_i\) (here \(\times\) denotes the categorical product).
The pullback can then be realised as the equaliser of the actors \(h_i\pi_i\), and hence we compute that.

Proposition~\ref{prop-productDisjUn} states that the product \(G_1\times G_2\) is realised as the disjoint union \(G_1\sqcup G_2\) with universal cone actors given by the relevant left-multiplications. 
The equaliser described in Subsection~\ref{subsec-EGALims} is the spatialisation of the subpseudogroup \(S=\{U\in\Bis (G_1\sqcup G_2): U\cdot_{\pi_1 h_1}H^{(0)}=U\cdot_{\pi_2 h_2}H^{(0)}\}\), which may be expressed as a quotient of the subgroupoid covered by \(S\).

Let \(K=\bigcup S\) be the open subgroupoid of \(G_1\sqcup G_2\) covered by \(S\).
Then \(K\) consists of elements \(\gamma\in G_i\) for which there are open bisections \(U_i\subseteq G_i\), \(V_j\subseteq G_j\) (where \(i\neq j\)) such that \(\gamma\in U_i\) and \(U_i\cdot_{h_i}H^{(0)}=V_j\cdot_{h_j} H^{(0)}\); the set \(\iota_i(U_i)\cup \iota_j(V_j)\subseteq G_1\sqcup G_2\) belongs to \(S\) (moreover every element of \(S\) arises this way).

We now form the unit space of pullback.
The anchor maps of the \(\pi_i h_i\) are given by composing the anchor map \(\rho_i\colon H\to G_i^{(0)}\) of \(h_i\) with the (\(\Top\)-canonical) inclusion \(\iota_i\colon G_i\hookrightarrow G_1\sqcup G_2\), and so the unit space of the pullback is the \(\Top\)-coequaliser of these maps.

For a unit \(t\in K^{(0)}\), the fibre over \([t]\in\coeq(\iota_1\rho_1,\iota_2\rho_2)\) in the pullback is represented by elements of the fibre over \(t\) in \(K\).
Within this fibre, two distinct arrows \(\gamma,\eta\in K_t\) are never identified, since no open bisection of \(K\) contains both as they have the same source.
In particular no element of \(S\) contains both.
This may be realised further as the quotient of \(K\) by the equivalence relation generated by declaring \(\gamma\sim\eta\) if \(s(\gamma)\sim s(\eta)\) and there is an element \(U\in S\) containing both \(\gamma\) and \(\eta\).

\subsection{A pullback of graph groupoids}

\subfile{graphPullback}

\subsection{Limits of groups}
The inclusion of the category \(\Grp\) of groups as a full subcategory of \(\EGA\) does not preserve all limits.
Indeed, Proposition~\ref{prop-productDisjUn} describes the product in \(\EGA\) as the disjoint union groupoid, which is seldom a group.
The underlying issue in this case, and for limits of diagrams of groups in \(\EGA\) more generally, is that a diagram in \(\EGA\) need not relate (or map) identities in distinct groups, and hence these identities may be separated by actors to other groupoids.

Recall that a category \(J\) is \emph{connected} if for any two objects \(a,b\in J\) there is a finite sequence of objects \(a=a_0,a_1,\dots,a_n=b\) together with arrows \(f_1,\dots,f_n\) in \(J\) where each \(f_i\) is either an arrow \(f_i\colon a_i\to a_{i+1}\) or \(f_i\colon a_{i+1}\to a_i\).
Such a sequence of arrows is called a \emph{zig-zag}.
We call a diagram \(F\colon J\to\Cc\) in a category \(\Cc\) \emph{connected} if \(J\) is connected.

\begin{proposition}\label{prop-GrpPresConnectedLims}
    The inclusion \(\iota\colon\Grp\to\EGA\) preserves connected limits.
    \begin{proof}
        Corollary~\ref{cor-EPresLims} implies that the unit space of a limit groupoid is given by the colimit in \(\Top\) of the diagram of the unit spaces and anchor maps.
        If the diagram is connected and all the unit spaces are singletons (which is exactly the case for a connected diagram of groups), the colimit of a connected diagram of singletons in \(\Top\) is a singleton, and so the limit of the diagram in \(\EGA\) has a singleton unit space i.e. is a group.
        Since the inclusion \(\iota\) is full and faithful, it reflects limits by \cite[Lemma~3.3.5]{Riehl_CategoryTheoryContext}, so if the limit of a diagram is a group then it was in fact the preserved limit from the category of groups.
    \end{proof}
\end{proposition}

Since \(\Grp\) has a terminal object, namely the trivial group, the limit of any diagram in \(\Grp\) may be realised as a connected limit by adjoining the terminal object.
Explicitly, if \(F\colon J\to\Grp\) is a small diagram of groups, define \(\tilde J\) as the category with objects \(J^\circ\sqcup\{\infty\}\) and arrows \(J^1\sqcup\{!\colon j\to\infty: j\in J^0\}\), and extend the diagram \(F\) to \(\tilde J\) by sending \(\infty \) to the trivial group \(\{e\}\) and each \(!\colon j\to \infty\) to the trivial map \(Fj\to\{e\}\).
Any cone over \(F\) extends naturally and uniquely to a cone over \(\tilde F\), so the corresponding limits are isomorphic in \(\Grp\).
Moreover, \(\tilde F\) is connected by construction, so this limit is preserved by the inclusion into \(\EGA\) by Proposition~\ref{prop-GrpPresConnectedLims}.

\begin{corollary}
    Let \(F\colon J\to\Grp\) be a small diagram.
    The limit of \(\iota F\colon J\to\EGA\) is the disjoint union of the limits of connected components of the diagram \(F\), in particular is the disjoint union of groups.
    \begin{proof}
        Express the shape category \(J\) as the disjoint union of its connected components \(J=\bigsqcup_{\alpha\in A}J_\alpha\).
        The diagram \(F\) then decomposes as the disjoint union of connected diagrams \(F_\alpha\colon J_\alpha\to\Grp\).
        A cone over \(F\) describes a family of cones over each \(F_\alpha\) by restriction, and conversely a family of cones over each \(F_\alpha\) with the same apex induces a cone over \(F\), yielding a natural isomorphism between cone functor \(\Cone(-,F)\) and the product \(\prod_{\alpha\in A}\Cone(-,F_\alpha)\).
        By the Yoneda Lemma, the limit of \(F\) is then the product of the limits of each \(F_\alpha\), and each \(F_\alpha\) is connected, so the limit of \(F_\alpha\) is preserved by \(\iota\) by Proposition~\ref{prop-GrpPresConnectedLims}.
        Finally, the product in \(\EGA\) is given by disjoint union by Proposition~\ref{prop-productDisjUn}, so the limit of \(\iota F\) is the disjoint union of the limits of each \(\iota F_\alpha\), each of which is a group by Proposition~\ref{prop-GrpPresConnectedLims}.
    \end{proof}
\end{corollary}

\subsection{A comment on colimits}

It is not clear whether colimits generally exist in  \(\Pseu\) or \(\EGA\).
Some specific cases are known, for example sequential inductive colimits are constructed in \cite{Tay_FunctGrpdCstarAlg}.
We also can quickly deduce that colimits of groups will behave well.

\begin{lemma}\label{lem-grpLeftAdj}
    The inclusions \(\Grp\hookrightarrow\Pseu\) adjoining a zero element and \(\Grp\hookrightarrow\EGA\) are left-adjoints.
    In particular, they preserve colimits.
    \begin{proof}
        For a group \(\Gamma\), denote the associated pseudogroup by \(\Gamma^0:=\Gamma\sqcup\{0\}\), where \(0\) is a zero element satisfying \(0t=t0=0\) for all \(t\in\Gamma^0\).
        For a pseudogroup \(S\), let \(S^\times\) be the set of invertible elements of \(S\).
        This is easily seen to be a group under the pseudogroup multiplication and inversion, and moreover any pseudogroup homomorphism \(\Gamma_0\to S\) restricts to a group homomorphism \(\Gamma\to S^\times\).
        We readily see that this restriction yields a family of bijections \(\Pseu(\Gamma_0,S)\cong \Grp(\Gamma,S^\times)\) natural in both \(\Gamma\) and \(S\).

        For groupoids, we observe that actors \(\Gamma\curvearrowright G\) from a group \(\Gamma\) to an \'etale groupoid \(G\) can be identified with group homomorphism \(\Gamma\to(\Bis G)^\times\).
        The natural bijections \(\EGA(\Gamma,G)\cong\Grp(\Gamma,(\Bis G)^\times)\) are then easily verified.
    \end{proof}
\end{lemma}

The above Lemma~\ref{lem-grpLeftAdj} makes use of the \emph{topological full group} of an \'etale groupoid \(G\); the group of global bisections of \(G\).
The topological full group has been studied extensively in the context of totally disconnected \'etale (or \emph{ample}) groupoids, even proving to be a complete invariant for such groupoids with no non-trivial closed invariant subspaces of units \cite{Matui_TFGsOneSidedShifts}.
Topological full groups have also been used to construct interesting examples of groups (see e.g. \cite{GardellaTanner_GeneraliseThompsonV}).

On the complete opposite end of the spectrum of groupoids, we also see that the coproduct of topological spaces viewed as \'etale groupoids corresponds to the product in \(\Top\).
Combining these two cases immediately complicates the issue, it is not obvious how to combine the limit structure of topological spaces with a fibrewise colimit or `free product'-like structure on fibres.

If a diagram in \(\Pseu\) admits a colimit, we begin to analyse the unit space of the corresponding spatialisation groupoid by applying the functor \(E\).

\begin{lemma}\label{lem-natEmbedd}
    Let \(S\) be a pseudogroup and let \(L\) be a frame.
    Any pseudogroup homomorphism \(\varphi\colon S\to\iota L\) is determined by the restriction to a frame homomorphism \(ES\to L\).
    Moreover, restriction describes natural embeddings \(\Pseu(S,\iota L)\hookrightarrow \Frame(ES,L)\).
    \begin{proof}
        Since \(\iota L\) is a pseudogroup consisting of idempotents, for any \(t\in S\) we have \(\varphi(t)=\varphi(t)^*\varphi(t)=\varphi(t^*t)\).
        Restriction is also natural for morphisms between frames, yielding a natural embedding \(\Pseu(S,\iota(-))\Rightarrow\Frame(ES,-)\).
    \end{proof}
\end{lemma}

If the colimit of such a diagram exists in \(\Pseu\), we may use the representable universal property of colimits to deduce that there is a natural frame epimorphism from the colimit of the idempotent frames to the idempotent frame of the colimit pseudogroup.

\begin{proposition}
    Let \(F\colon J\to\Pseu_\Spat\) be a diagram which admits a colimit in \(\Pseu_\Spat\).
    Then there is a frame epimorphism \(\colim EF\to E\colim F\).
    Dually, there is a continuous injection 
    \[(\colim\Sigma F)^{(0)}\hookrightarrow \lim\limits_{J^\op,\Top}(\Sigma Fj)^{(0)},\] 
    where the limit on the right side is taken in the category of topological spaces of the diagram consisting of anchor maps of the morphisms in the diagram \(F\).
    \begin{proof}
        Let \(F\colon J\to\Pseu_\Spat\) be a diagram admitting a colimit \(\colim F\).
        Using Proposition~\ref{prop-frmCorefPseu}, Lemma~\ref{lem-natEmbedd} and the representable universal property of colimits (see \cite[Theorem~3.4.7]{Riehl_CategoryTheoryContext})
        \begin{align*}
            \Frame(E\colim F,-)&\cong \Pseu(\colim F,\iota(-))\\
            &\cong \lim\limits_{J^\op}\Pseu(Fj,\iota(-))\\
            &\hookrightarrow \lim\limits_{J^\op}\Frame(EFj,-)\\
            &\cong \Frame(\colim EF,-).
        \end{align*}
        By the Yoneda Lemma, there is a unique epimorphism \(\alpha\colon \colim EF\to E\colim F\) such that precomposition with \(\alpha\) induces this natural transformation.
        After applying the spatialisation functor, the unit space of \(\Sigma\colim F\cong\colim\Sigma F\) is the spectrum of \(\colim EF\), since \(\Sigma\) is a left-adjoint so preserves colimits.
        Hence we get the desired monomorphism.
    \end{proof}
\end{proposition}

%% file: sections/graphPullback.tex
We now give an explicit example of a pullback computation for groupoids arising from directed graphs.

Given a directed graph \(E=(E^0,E^1,r,s)\) with vertex set \(E^0\), edge set \(E^1\), and range and source maps \(r,s\colon E^1\to E^0\), Paterson \cite{Paterson_GraphGrpd} defines the graph groupoid \(G_E\).
We recall an equivalent construction here.

A \emph{path} in \(E\) is a finite or infinite string \(e_1\cdots e_n\) or \(e_1\cdots\) of edges \(e_i\in E^1\) such that \(s(e_i)=r(e_{i+1})\) for \(i\geq 1\leq n\) (or simply \(i\geq 1\) if the path is infinite).
We take the convention that paths read right-to-left to match with our convention on composition in groupoids and pseudogroups.
The \emph{length} of a path \(x=e_1\cdots e_n\) is \(n\), and is denoted \(|x|\).
We also allow vertices as paths of length zero.
We write \(E^*\) for the set of finite paths in \(E\), and \(E^\infty\) for the set of infinite paths.
We extend the source map on \(E\) to \(E^*\) by setting \(s(e_1\cdots e_n)=s(e_n)\), and we extend the range map to \(E^*\cup E^\infty\) by setting \(r(e_1\cdots)=r(e_1)\).
A vertex \(v\in E^0\) is a \emph{source} if \(r^{-1}(v)=\emptyset\), that is, \(v\) has no incoming edges.
The unit space of \(G_E\) is the \emph{path space}
\[\partial E=E^\infty\cup \{\lambda\in E^*: s(\lambda)\text{ is a source}\},\]
with topology given by basic open sets \(Z_x=\{\lambda\in \partial E: \lambda=x\lambda'\text{ for some}\lambda'\in\partial E\}\) indexed over finite paths \(x\in E^*\).
Here \(x\lambda'\) denotes the concatination of the two paths \(x\) and \(\lambda'\), which is defined if and only if \(s(x)=r(\lambda)\).

Let \(G_E\) be the subset of \(\partial E\times\ZZ\times\partial E\) consisting of triples \((\lambda,n-m,\mu)\), where \(\lambda\) and \(\mu\) are paths which coincide after the deleting the left-most \(n\) edges from \(\lambda\) and \(m\) edges from \(\mu\).
If \(n\) is larger than the length of \(\lambda\) (in case \(\lambda\) is finite), deleting the left-most \(n\) edges results in the zero-length path consisting of the source of \(\lambda\).

We equip \(G_E\) with the subspace topology, which we note admits a base consisting of sets of the form 
\begin{equation}\label{eqn-basicBisections}
    Z^E(x,y):=\{(x\lambda,|x|-|y|,y\lambda): \lambda\in \partial E, s(x)=s(y)=r(\lambda)\}
\end{equation}
for finite paths \(x,y\in E^*\).
We may write \(Z(x,y)\) for \(Z^E(x,y)\) if the graph \(E\) is unambiguous.
This turns \(G_E\) into a locally compact Hausdorff space.
In particular, it is sober.

We identify \(\partial E\) with its image in \(G_E\) under the map \(\lambda\mapsto(\lambda,0,\lambda)\).
Define range and source maps on \(G_E\) by \(r(\lambda,n,\mu)=\lambda\) and \(s(\lambda,n,\mu)=\mu\), and the composition \(G_E\baltimes{s}{r} G_E\to G_E\) by 
\[(\lambda,n,\mu)(\mu,m,\nu)=(\lambda,n+m,\nu).\]
These operations turn \(G_E\) into a topological groupoid with unit space given \(\partial E\).
The basic open sets \(Z(x,y)\) in Equation~\eqref{eqn-basicBisections} form open bisections of \(G_E\), implying that \(G_E\) is \'etale.

The bisection pseudogroup associated to the groupoid \(G_E\) is generated by bisections of the form \(Z(e,s(e))\) for edges \(e\in E^1\).
These bisections encode the dynamical structure of the graph \(E\), leading to the folloiwng definition.

\begin{definition}[{\cite[Definition~1.1]{deCastroMeyer_GraphMorphismsGrpdActors}}]
    Let \(E\) be a directed graph and let \(H\) be an \'etale groupoid.
    A \emph{dynamical Cuntz-Krieger \(E\)-family} in \(H\) is a family of compact open subsets \(\Omega=(\Omega_v)_{v\in E^{(0)}}\) and a family of compact open bisections \(T=(T_e)_{e\in E^1}\) in \(H\) satisfying
    \begin{enumerate}
        \item \(\Omega_v\cap\Omega_w=\emptyset\) for distinct vertices \(v,w\in E^0\),
        \item \(s(T_e)\subseteq\Omega_{s(e)}\) and \(r(T_e)\subseteq \Omega_{r(e)}\) for all edges \(e\in E^1\),
        \item \(r(T_e)\cap r(T_f)=\emptyset\) for distinct edges \(e,f\in E^1\),
        \item \(\Omega_v=\bigsqcup_{r(e)=v}r(T_e)\) for any vertex \(v\in E^0\) with \(r^{-1}(v)\) finite and non-empty.
    \end{enumerate}
\end{definition}

De Castro and Meyer show that there are natural bijections between actors \(G_E\curvearrowright H\) and Cuntz-Krieger \(E\)-families in \(H\).

\begin{theorem}[{\cite[Theorem~4.3, Corollary~4.6]{deCastroMeyer_GraphMorphismsGrpdActors}}]\label{thm-graphThm}
    Any Cuntz-Krieger \(E\)-family \((\Omega,T)\) in \(H\) induces an actor \(G_E\curvearrowright H\) satisfying \(Z_v\cdot H^{(0)}=\Omega_v\) and \(Z(e,s(e))\cdot H^{(0)}=T_e\).
    The functions sending Cuntz-Krieger \(E\)-families to the associated actors are natural bijections between the sets of Cuntz-Krieger \(E\)-families and the actors with domain \(G_E\).
\end{theorem}

We shall use Theorem~\ref{thm-spatSobComplete} and Theorem~\ref{thm-graphThm} to describe a particular pullback of graph groupoids.

Consider the graph \(E\) with vertices \(E^0=\{v_1,v_2,w\}\), edges \(E^1=\{e_1,e_2,\sigma\}\), and range and source maps given by \(r(e_i)=r(\sigma)=w\), \(s(e_i)=v_i\), and \(s(\sigma)=w\):

\[\begin{tikzcd}
	& w \\
	{v_1} && {v_2}
	\arrow["\sigma", from=1-2, to=1-2, loop, in=55, out=125, distance=10mm]
	\arrow["{e_1}", from=2-1, to=1-2]
	\arrow["{e_2}"', from=2-3, to=1-2]
\end{tikzcd}\]

Set \(E_1\) to be the full subgraph on the vertices \(\{v_1,w\}\) and \(E_2\) to be the full subgraph on vertices \(\{v_2,w\}\).
Set \(F\) to be the full subgraph on the vertex \(w\) which consists of single point and loop.

\[\begin{tikzcd}
	{E_1:} && w && {E_2:} & w \\
	& {v_1} &&&&& {v_2} \\
	\\
	&& {F:} & w
	\arrow["\sigma", from=1-3, to=1-3, loop, in=55, out=125, distance=10mm]
	\arrow["\sigma", from=1-6, to=1-6, loop, in=55, out=125, distance=10mm]
	\arrow["{e_1}", from=2-2, to=1-3]
	\arrow["{e_2}"', from=2-7, to=1-6]
	\arrow["\sigma", from=4-4, to=4-4, loop, in=55, out=125, distance=10mm]
\end{tikzcd}\]

The choices \(\Omega^i_{w}=Z_w\), \(T^i_\sigma=Z^F(\sigma,s(\sigma))\), and \(T^i_{e_i}=\Omega_{v_i}=\emptyset\) define Cuntz-Krieger \(E_i\)-families \((\Omega^i,T^i)\) in \(G_F\) inducing actors \(G_{E_1},G_{E_2}\curvearrowright G_F\) by Theorem~\ref{thm-graphThm}.
We aim to compute the limit of this pullback diagram using the machinery we developed in this article.
Since graph groupoids are sober, it suffices to consider the corresponding diagram of pseudogroups 
\begin{equation*}\label{diag-pullbackPseu}
    \begin{tikzcd}
        & {\Bis(G_{E_1})} \\
        {\Bis(G_{E_2})} & {\Bis(G_{F})}
        \arrow["{\psi_1}", from=1-2, to=2-2]
        \arrow["{\psi_2}"', from=2-1, to=2-2]
    \end{tikzcd}
\end{equation*}
where \(\psi_i\) is the pseudogroup homomorphism induced by the Cuntz-Krieger \(E_i\)-families in \(G_F\).
The \(\psi_i\) are the unique pseudogroup homomorphisms satisfying \(\psi_i(Z^{E_i}(e,s(e)))=T^i_e\) for each edge \(e\in E^1_i\).
Using Theorem~\ref{thm-pseuComplete}, we can realise the limit of \eqref{diag-pullbackPseu} as the set
\[\{(U_1,U_2)\in\Bis(G_{E_1})\times\Bis(G_{E_2}):\psi_1(U_1)=\psi_2(U_2)\}\]
equipped with pointwise operations.

\begin{lemma}\label{lem-removeInfSentZero}
    We have \(\psi_i(\partial E_i\setminus\{\sigma^\infty\})=\emptyset\) for \(i=1,2\).
    \begin{proof}
        We may express \(\partial E_i\setminus\{\sigma^\infty\}\) as the union of the sources and ranges of the bisections \(Z^{E_i}(\sigma,s(\sigma))^n Z^{E_i}(e_i,s(e_i))\) over all \(n\in\NN_0\).
        The pseudogroup morphisms \(\psi_i\) send \(Z^{E_i}(e_i,s(e_i))\) to the empty set, and each \(\psi_i\) preserves products, yielding the lemma. 
    \end{proof}
\end{lemma}

\begin{lemma}\label{lem-psiSendsSingleton}
    Fix \(U_i\in\Bis(G_{E_i})\).
    The bisection \(\psi_i(U_i)\in\Bis(G_F)\) is either the singleton containing the (necessarily unique) element of \(U_i\) with source and range equal to \(\sigma^\infty\), or otherwise is empty if \(U_i\) contains no such element.
    \begin{proof}
        If \(U_i\) does not contain a point with both source and range given by \(\sigma^\infty\) then we have either \(U_i=U_i(\partial E_i\setminus\{\sigma^\infty\})\) or \(U_i=(\partial E_i\setminus\{\sigma^\infty\})U_i\).
        Lemma~\ref{lem-removeInfSentZero} implies \(\psi_i(U_i)=\emptyset\) in both cases.

        Now suppose \(U_i\) contains an element with source and range given by \(\sigma^\infty\).
        Necessarily this element has the form \((\sigma^\infty,n,\sigma^\infty)\) for some \(n\in\ZZ\).
        By replacing \(U_i\) with its adjoint if necessariy, we may assume \(n\geq 0\).
        We can find a basic neighbourhood of the point \((\sigma^\infty,n,\sigma^\infty)\) of the form \(Z(x,y)\) for some \(x,y\in E_i^*\) contained in \(U_i\), and we deduce necessarily that \(x=\sigma^{k+n}\) and \(y=\sigma^k\) for some \(k\in\NN_0\).
        Lemma~\ref{lem-removeInfSentZero} implies that
        \[\psi_i(U_i)=\psi_i(U_i(\partial E_i\setminus\{\sigma^\infty\})\cup Z(x,y))=\psi(Z(x,y)).\]
        Observe that \(Z(\sigma^{n+k},\sigma^k)\) can be expressed as the product \(Z(\sigma,s(\sigma))^{n+k}Z(\sigma,s(\sigma))^{k*}\) and \(\psi_i\) sends \(Z(\sigma,s(\sigma))\) to the singleton \(\{(\sigma^\infty,1,\sigma^\infty)\}\subseteq G_F\), so we have
        \[\psi_i(U_i)=\psi(Z(\sigma,s(\sigma)))^{n+k}\psi(Z(\sigma,s(\sigma)))^{k*}=\{(\sigma^\infty,n,\sigma^\infty)\}\]
        as required.
    \end{proof}
\end{lemma}

\begin{lemma}\label{lem-generatorsGenerate}
    The pullback \(\Bis(G_{E_1})\times_{\Bis(G_F)}\Bis(G_{E_2})\) of the diagram \eqref{diag-pullbackPseu} is generated by the elements 
    \begin{itemize}
        \item \((Z^{E_1}(e_1,s(e_1)),\emptyset)\), 
        \item \((\emptyset,Z^{E_2}(e_2,s(e_2)))\), 
        \item \((Z^{E_1}(\sigma,s(\sigma)),Z^{E_2}(\sigma,s(\sigma)))\).
    \end{itemize}
    \begin{proof}
        Write \(S:=\Bis(G_{E_1})\times_{\Bis(G_F)}\Bis(G_{E_2})\).
        We first note that each of the proposed generators lies in \(S\) since \(\psi_i(Z^{E_i}(e_i,s(e_i)))=T^i_{e_i}=\emptyset\) and \(\psi_1(Z^{E_1}(\sigma,s(\sigma)))=Z^F(\sigma,s(\sigma))=\psi_2(Z^{E_2}(\sigma,s(\sigma)))\).
        Let \(S'\subseteq S\) be the subpseudogroup generated by the proposed generators.

        As in Lemma~\ref{lem-removeInfSentZero} we may express \(\partial E_i\setminus\{\sigma^\infty\}\) as the union of the sources and ranges of the bisections \(Z^{E_i}(\sigma,s(\sigma))^n Z^{E_i}(e_i,s(e_i))\) over all \(n\in\NN_0\).
        All of these components are sent to the empty set under \(\psi_i\), since by Lemma~\ref{lem-removeInfSentZero} we have\(\psi_i(\partial E_i\setminus\{\sigma^\infty\})=\emptyset=\psi_j(\emptyset)\) for \(j\neq i\).
        Hence \((\partial E_1\setminus\{\sigma^\infty\},\emptyset)\) and \((\emptyset,\partial E_2\setminus\{\sigma^\infty\})\) are both elements of the \(S\) which we can generate from the generators of \(S'\).
        Write \(G_{E_i}':=G_{E_i}\cdot \partial E_i\setminus\{\sigma^\infty\}\) and note that this is an open subgroupoid of \(G_{E_i}\) for \(i=1,2\).
        This implies that \(\Bis(G_{E_i}')\) is an inverse subsemigroup of \(\Bis(G_{E_i})\), and further \(\psi_i\) sends all bisections of \(G_{E_i}'\) to the emptyset, so \((U,\emptyset)\) belongs to \(S\) for all \(U\in\Bis(G_{E_i}')\). 
        The bisections \(Z^{E_i}(e_i,s(e_i))\cdot(\partial E_i\setminus \{\sigma^\infty\})\) and \(Z^{E_i}(\sigma,s(\sigma))\cdot (\partial E_i\setminus\{\sigma^\infty\})\) then generate \(\Bis(G_{E_i}')\), and so we see that the product \(\Bis(G_{E_1}')\times\Bis(G_{E_2}')\) embeds in the obvious way into \(S\) as an inverse subsemigroup, and moreover lies in \(S'\).

        Fix \((U_1,U_2)\in S\).
        We note that if \(\psi_1(U_1)=\emptyset=\psi_2(U_2)\) then \(U_i\subseteq G_{E_i}'\), and so \((U_1,U_2)\) belongs to \(S'\).
        We now consider the case where \(\psi_1(U_1)=\psi_2(U_2)\neq\emptyset\).
        Lemma~\ref{lem-psiSendsSingleton} then implies that there is \(n\in\ZZ\) such that \(\gamma:=(\sigma^\infty,n,\sigma^\infty)\in U_1\cap U_2\).
        We can then find \(x,y\in E_1^*\) and \(z,w\in E_2^*\) such that \(\gamma\in Z^{E_1}(x,y)\subseteq U_1\) and \(\gamma\in Z^{E_2}(z,w)\in U_2\).
        Since elements of e.g. \(Z^{E_1}(x,y)\) are of the form \(x\lambda,|x|-|y|,y\lambda\) for some \(\lambda\in \partial E_1\) and \(\gamma=(\sigma^\infty,n,\sigma^\infty)\) belongs to this set, it follows that \(x=\sigma^{|x|}\), \(y=\sigma^{|y|}\), and similarly for \(z\) and \(w\).
        Moreover we have \(|x|-|y|=n=|z|-|w|\).
        By interchanging \((U_1,U_2)\) for \((U_1,U_2)^*=(U_1^*,U_2^*)\) if necessary, it suffices to consider the case where \(n\geq 0\), hence \(k\geq n\).
        We then have \(Z^{E_1}(\sigma^k,\sigma^{n-k})\subseteq Z^{E_1}(x,y)\) and \(Z^{E_2}(\sigma^k,\sigma^{n-k})\subseteq Z^{E_2}(w,z)\), and \(\gamma\) belongs to both.
        Moreover \((Z^{E_1}(\sigma^k,\sigma^{n-k}),Z^{E_1}(\sigma^k,\sigma^{n-k}))\) can be expressed as products of \((Z^{E_1}(\sigma,s(\sigma)),Z^{E_2}(\sigma,s(\sigma)))\) and its inverse, so the pair \((Z^{E_1}(\sigma^k,\sigma^{n-k}),Z^{E_1}(\sigma^k,\sigma^{n-k}))\) belongs to \(S'\).
        Finally, we note that \(\gamma\) is the only element of \(U_1\) and \(U_2\) which does not belong to \(U_i\cdot (\partial E_i\setminus\{\sigma^\infty\})\).
        Hence we have
        \begin{align*}
            (U_1,U_2)=&(U_1\cdot(\partial E_1\setminus\{\sigma^\infty\}),U_2\cdot(\partial E_2\setminus\{\sigma^\infty\}))\\
            &\vee(Z^{E_1}(\sigma^k,\sigma^{n-k}),Z^{E_2}(\sigma^k,\sigma^{n-k}))
        \end{align*}
        so \((U_1,U_2)\) is contained in \(S'\) since both of these components are.
    \end{proof}
\end{lemma}

\begin{theorem}\label{thm-graphPullback}
    The pullback of the actors \(G_{E_i}\curvearrowright G_F\) for \(i=1,2\) induced by the above Cuntz-Krieger \(E_i\)-families in \(F\) is the graph groupoid \(G_E\).
    The universal cone actors \(G_E\curvearrowright G_{E_i}\) are induced by dynamical Cuntz-Krieger \(E\)-families \((\Omega^i,T^i)\) in \(E_i\) given by
    \begin{align*}
        \Omega^i_w&:=Z_w,& T^i_{\sigma}&:=Z(\sigma,s(\sigma)),\\
        \Omega^i_{v_j}&:=\begin{cases}
            Z_{v_i},& i=j,\\
            \emptyset,&\text{otherwise,}
        \end{cases}&T^i_{e_j}&:=\begin{cases}
            Z(e_i,s(e_i)),&i=j,\\
            \emptyset,&\text{otherwise.}
        \end{cases}
    \end{align*}
    \begin{proof}\
        Denote the pullback \(\Bis(G_{E_1})\times_{\Bis(G_{F})}\Bis(G_{E_2})\) by \(S\).
        Let \(\varphi_i\colon \Bis(G_E)\to\Bis(G_{E_i})\) be the pseudogroup homomorphisms induced by the dynamical Cuntz-Krieger \(E\)-families \((\Omega^i,T^i)\).
        We readily see that since the compositions \(\psi_1\varphi_1\) and \(\psi_2\varphi_2\) coincide by evaluating on the canonical generators of \(\Bis(G_E)\), hence there is a unique pseudogroup homomorphism \(\varphi\colon \Bis(G_E)\to S\) satisfying \(\pi_i\varphi=\varphi_i\) for \(i=1,2\), where \(\pi_i\colon S\to \Bis(G_{E_i})\) is the canonical projection.

        We denote the powerset of \(G_E\) by \(\Pp(G_E)\).
        By viewing \(E_1\) and \(E_2\) directly as subgraphs of \(E\), we obtain inclusion maps \(\iota_i\colon\Bis(G_{E_i})\to\Pp(G_E)\) from the bisection pseudogroups to the powerset of the graph groupoid \(G_E\) sending any bisection of \(G_{E_i}\) to the set itself, viewed as a subset of \(G_E\).
        We may similarly consider \(\Bis(G_E)\) as a subset of \(\Pp(G_E)\).
        Using the characterisation of \(S\) as the set
        \[S:=\{(U_1,U_2)\in\Bis(G_{E_1})\times\Bis(G_{E_2}): \psi_1(U_1)=\psi_2(U_2)\},\]
        we claim the function \(\alpha\colon S\to\Pp(G_E)\) given by
        \[\alpha(U_1,U_2)=\iota_1(U_1)\cup\iota_2(U_2)\]
        takes values in \(\Bis(G_E)\) and defines a pseudogroup homomorphism inverse to \(\varphi\).

        Since the source and range maps on \(G_E\) restrict to the source and range maps on each \(G_{E_i}\), we note that the range and source maps are injective on each \(\iota_i(U_i)\) since the \(U_i\) are bisections in \(G_{E_i}\).
        Suppose there are \(\gamma_1\in \iota_1(U_1)\) and \(\gamma_2\in\iota_2(U_2)\) have the same source.
        This implies that \(s(\gamma_1)=s(\gamma_2)\) belongs to the intersection \(\partial E_1\cap \partial E_2=\{\sigma^\infty\}\), so \(\gamma_1\) and \(\gamma_2\) have source \(\sigma^\infty\).
        Since \(\sigma^\infty\) is an invariant point in both graph groupoids \(G_{E_i}\), we deduce that \(r(\gamma_1)=\sigma^\infty=r(\gamma_2)\), so \(\gamma_i\in U_i\) is the unique element of \(U_i\) which has source and range \(\sigma^\infty\).
        Applying Lemma~\ref{lem-psiSendsSingleton} then yields \(\{\gamma_1\}=\psi_1(U_1)=\psi_2(U_2)=\{\gamma_2\}\), implying \(\gamma_1=\gamma_2\).
        Hence \(s\) is injective on \(\iota_1(U_1)\cup\iota_2(U_2)\), and a symmetric argument shows that \(r\) is as well.

        To see that \(\alpha(U_1,U_2)\) is open, fix a point \(\gamma\in\alpha(U_1,U_2)\).
        Since the intersection of the groupoids \(G_{E_i}\) (viewed as subgroupoids of \(G_E\)) contains exactly the fibre over \(\sigma^\infty\), if \(s(\gamma)\neq\sigma^\infty\) then \(\gamma\) belongs to \(U_i(\partial E_i\setminus\{E_i\})\) for one \(i=1,2\).
        Since \(\partial E_i\setminus\{\sigma^\infty\}\) is open in \(\partial E\), we see that \(\iota_i(U_i(\partial E_i\setminus\{\sigma^\infty\}))\) is an open neighbourhood of \(\gamma\) contained in \(\alpha(U_1,U_2)\).
        Otherwise if \(s(\gamma)=\sigma^\infty\) then we can apply Lemma~\ref{lem-psiSendsSingleton} to see \(\gamma\in U_1\cap U_2\) and as before can find a neighbourhood \(Z^E(\sigma^k,\sigma^\ell)\) of \(\gamma\) contained in both \(U_1\) and \(U_2\).
        Hence \(\alpha(U_1,U_2)\) is an open bisection, so \(\alpha\) takes values in \(\Bis(G_E)\).

        We note that the maps \(\iota_1\) and \(\iota_2\) are \({}^*\)-preserving and the \({}^*\)-operation preserves unions, so \(\alpha\) is \({}^*\)-preserving.
        To see \(\alpha\) is multiplicative, fix \((U_1,U_2),(V_1,V_2)\in S\) and observe that pointwise product \(\iota_1(U_1)\iota_2(V_2)\) os either empty or contains only the product of the unique elements of \(U_1\) and \(V_2\) with source and range equal to \(\sigma^\infty\).
        Lemma~\ref{lem-psiSendsSingleton} then implies that \(\iota_1(U_1)\iota_1(V_1)\) also contains this singleton since \(\psi_1(V_1)=\psi_2(V_2)\), so \(\iota_1(U_1)\iota_2(V_2)\subseteq \iota_1(U_1)\iota_1(V_1)\).
        Symmetrically we have \(\iota_2(U_2)\iota_1(V_1)\subseteq\iota_2(U_2)\iota_2(V_2)\).
        We also note that each \(\iota_i\) is multiplicative, so we have
        \begin{align*}
            \alpha(U_1,U_2)\alpha(V_1,V_2)&=(\iota_1(U_1)\cup\iota_2(U_2))(\iota_1(V_1)\cup\iota_2(V_2))\\
            &=\iota_1(U_1)\iota_2(U_2)\cup\iota_1(U_1)\iota_2(U_2)\cup\iota_2(U_2)\iota_1(V_1)\cup\iota_2(U_2)\iota_2(V_2)\\
            &=\iota_1(U_1)\iota_1(V_1)\cup\iota_2(U_2)\iota_2(V_2)\\
            &=\iota_1(U_1V_1)\cup\iota_2(U_2V_2)\\
            &=\alpha(U_1V_1,U_2V_2),
        \end{align*}
        so \(\alpha\) is multiplicative.
        We also note that \(\alpha\) preserves (compatible) joins since the inclusion maps \(\iota_i\) preserve unions and \(\alpha\) is defined as a union of these.
        Hence \(\alpha\) is a pseudogroup morphism.

        Finally, we demonstrate that \(\alpha\) is the inverse to \(\varphi\).
        The pseudogroup \(\Bis(G_E)\) is generated by the bisections \(Z^E(e_i,s(e_i))\) and \(Z^E(\sigma,s(\sigma))\), while Lemma~\ref{lem-generatorsGenerate} shows \(S\) is generated by \((\emptyset,Z^{E_i}(e_i,s(e_i)))\), and \((Z^{E_1}(\sigma,s(\sigma)),Z^{E_2}(\sigma,s(\sigma)))\).
        By construction of \(\alpha\) and since \(\varphi\) is induced by the Cuntz-Krieger \(E\)-families \((\Omega^i,T^i)\), we see that \(\varphi\) and \(\alpha\) restrict to mutually inverse bijections on these sets of generators.
        Hence they form mutually inverse isomorphisms.
        Corollary~\ref{cor-catEquivSobSpat} then implies that \(G_E\) is the corresponding pullback groupoid in \(\EGA\).
    \end{proof}
\end{theorem}

In \cite{HajacReznikoffTobolski_PullbacksGraphCStarAlgs}, Hajac, Reznikoff, and Tobolski specify criteria for when a pushout of directed graphs lifts to a pullback of the corresponding \(C^*\)-algebras.
They give the above graphs \(F\subseteq E_1,E_2\subseteq E\) as one such example.
Theorem~\ref{thm-graphPullback} implies that the graph \(E\) also induces the corresponding pullback of groupoids and actors.

%% file: appendices/catTheory.tex
For a general reference on category theory we refer the reader to \cite{MacLane_WorkingMathematician} or \cite{Riehl_CategoryTheoryContext}.

\subsection{Elementary definitions, the Yoneda Lemma}

\begin{definition}
    A \emph{category} \(\Cc\) consists of
    \begin{itemize}
        \item a class \(\Cc^0\) of \emph{objects},
        \item a class \(\Cc^1\) of \emph{arrows} or \emph{morphisms} between objects in \(\Cc\).
        Each arrow \(f\) has a unique domain \(c\) and codomain \(d\), denoted \(f\colon c\to d\), 
        \item \emph{composition}: if \(f\colon c\to d\) and \(g\colon d\to e\) are arrows in \(\Cc\) then there is a composition arrow \(gf\colon c\to e\),
    \end{itemize}
    such that 
    \begin{itemize}
        \item for each object \(c\in\Cc\) there is an \emph{identity arrow} \(1_c\colon c\to c\) satisfying \(f1_c=f\) and \(1_cg=g\) for all \(f\colon c\to d\) and \(g\colon d\to c\),
        \item composition is associative: \((hg)f=h(gf)\) whenever \(f\colon c\to d\), \(g\colon d\to e\), and \(h\colon e\to b\) are composable arrows in \(\Cc\).
    \end{itemize}
    An arrow \(f\colon c\to d\) is an \emph{isomorphism} if there is an arrow \(g\colon d\to c\) satisfying \(gf=1_c\) and \(fg=1_d\).

    A category \(G\) is a \emph{groupoid} if the morphism class \(G=G^1\) is a set and every morphism is an isomorphism.
\end{definition}

Every category \(\Cc\) has an \emph{opposite} category \(\Cc^\op\) with the same objects in \(\Cc\) and where every morphism \(f\colon c\to d\) in \(\Cc\) describes a morphism \(f^\op\colon d\to c\).
The composition in \(\Cc^\op\) is given by \(g^\op f^\op=(fg)^\op\).

\begin{definition}
    A (\emph{covariant}) \emph{functor} \(F\colon C\to D\) between two categories \(C\) and \(D\) consists of functions \(\Cc^0\to\Dd^0\) and \(\Cc^1\to \Dd^1\) (both denoted by \(F\)) satisfying 
    \begin{itemize}
        \item \(Ff\colon Fc\to Fd\) for all morphisms \(f\colon c\to d\) in \(\Cc\),
        \item \(F1_c=1_{Fc}\) for all objects \(c\) in \(\Cc\).
    \end{itemize}
    A \emph{contravariant functor} from \(\Cc\) to \(\Dd\) is a covariant functor \(\Cc^\op\to\Dd\).
\end{definition}

Functors \(\Cc\to\Dd\) and \(\Dd\to\Ee\) may be composed to get a functor \(\Cc\to\Ee\).
Each category \(\Cc\) admits an identity functor \(1_\Cc\) which consists of the identity maps on \(\Cc^0\) and \(\Cc^1\).
This functor fills the role of an identity with repsect to functor composition.

\begin{definition}
    Given functors \(F,G\colon \Cc\to\Dd\) between categories, a \emph{natural transformation} \(\alpha\colon F\Rightarrow G\) from \(F\) to \(G\) is a family \((\alpha_c)_{c\in\Cc^0}\) of arrows \(\alpha_c\colon Fc\to Gc\) in \(\Dd\), called \emph{components}, such that for each arrow \(f\colon c\to d\) in \(\Cc\) we have \(\alpha_d Ff=Gf\alpha_c\).
    Expressed diagrammatically, the diagram
    \[\begin{tikzcd}
        Fc & Gc \\
        Fd & Gd
        \arrow["{\alpha_c}", from=1-1, to=1-2]
        \arrow["Ff"', from=1-1, to=2-1]
        \arrow["Gf", from=1-2, to=2-2]
        \arrow["{\alpha_d}"', from=2-1, to=2-2]
    \end{tikzcd}\]
    commutes.
    A natural transformation \(\alpha\) is a \emph{natural isomorphism} if each component \(\alpha_c\) is an isomorphism.
    For functors \(F,G\colon \Cc\to\Dd\), we write \(\Nat(F,G)\) for the class of natural transformations \(F\Rightarrow G\).
\end{definition}

Natural transformations \(\alpha\colon F\Rightarrow G\) and \(\beta\colon G\Rightarrow H\) may be composed component-wise yielding a natural transformation \(\beta\alpha\colon F\Rightarrow H\) with components \((\beta\alpha)_c:=\beta_c\alpha_c\).

Given a category \(\Cc\) and objects \(c,d\in\Cc^0\), we write \(\Cc(c,d)\) for the class of morphisms \(f\colon c\to d\).
We say \(\Cc\) is \emph{locally small} if \(\Cc(c,d)\) is a set for all objects \(c,d\in\Cc^0\).
We denote the category of sets and functions by \(\Set\).
If \(\Cc\) is locally small, each object \(c\in\Cc^0\) defines two functors \(\Cc(c,-)\colon \Cc\to\Set\) and \(\Cc(-,c)\colon \Cc^\op\to\Set\).
On objects, \(\Cc(c,-)d:=\Cc(c,d)\) and \(\Cc(-,c)d=\Cc(d,c)\), and an arrow \(f\colon d\to e\) defines functions \(f_*\colon\Cc(c,d)\to\Cc(c,e)\) and \(f^*\colon\Cc(e,c)\to\Cc(d,c)\) by post and precomposition respectively.
A \emph{representation} of a covariant functor \(F\colon\Cc\to\Set\) is a natural isomorphism \(\Cc(c,-)\Rightarrow F\).
The functor \(F\) is \emph{representable} if it admits a representation \(\Cc(c,-)\Rightarrow F\), in this case the object \(c\) is called a \emph{representative} of \(F\).
The Yoneda Lemma helps to describe natural transformations from represented functors to other functors.

\begin{theorem}[Yoneda Lemma {\cite[III.2(Lemma)]{MacLane_WorkingMathematician}}]
    Let \(\Cc\) be a locally small category and let \(F\colon\Cc\to\Set\) be a functor.
    The function \(y\colon \Nat(\Cc(c,-),F)\to Fc\) sending a natural transformation \(\alpha\colon \Cc(c,-)\Rightarrow F\) to \(y(\alpha):=\alpha_c(1_c)\) is a bijection for all \(c\in\Cc^0\).
    Moreover, these bijections are natural in both \(c\) and \(F\).
\end{theorem}

An immediate consequence of the Yoneda Lemma is that all natural transformations between represented functors are implemented by morphisms in the category.

\subsection{Limits}

A \emph{diagram} of shape \(J\) in a category \(\Cc\) is a functor \(F\colon J\to \Cc\).
A diagram is \emph{small} if the shape category \(J\) is.

While there is formally no difference between a diagram and a functor, the typical use of each name is different.
Diagrams are often used in the contexts of limits.

\begin{definition}
    Let \(F\colon J\to \Cc\) be a diagram and let \(c\in\Cc\) be an object.
    A \emph{cone over \(F\)} with apex \(c\) is a family \(\lambda=(\lambda_j)_{j\in J}\) of morphisms \(\lambda_j\colon c\to Fj\) in \(\Cc\) for each \(j\in J\), called \emph{legs} or \emph{components}, satisfying \(\lambda_k=(Ff)\lambda_j\) for all arrows \(f\colon j\to k\) in \(J\).
    \[\begin{tikzcd}
        & c \\
        Fj && Fk
        \arrow["{\lambda_j}"', from=1-2, to=2-1]
        \arrow["{\lambda_k}", from=1-2, to=2-3]
        \arrow["Ff"', from=2-1, to=2-3]
    \end{tikzcd}\]
    A cone \(\lambda\) over \(F\) with apex \(c\) is a \emph{limit cone} if for every cone \(\mu\) over \(F\) with apex \(d\), there is a unique arrow \(\tilde\mu\colon d\to c\) satisfying \(\mu_j=\lambda_j\mu\) for all \(j\in J\).
    \[\begin{tikzcd}
        & d \\
        & c \\
        Fj && Fk
        \arrow["{\exists ! \tilde\mu}"', dashed, from=1-2, to=2-2]
        \arrow["{\mu_j}"', bend right, from=1-2, to=3-1]
        \arrow["{\mu_k}", bend left, from=1-2, to=3-3]
        \arrow["{\lambda_j}"', from=2-2, to=3-1]
        \arrow["{\lambda_k}", from=2-2, to=3-3]
        \arrow["Ff"', from=3-1, to=3-3]
    \end{tikzcd}\]
    In this case we call \(c\) the \emph{limit} of \(F\).
\end{definition}

Cones under diagrams and colimits are defined dually.
A cone under a diagram \(F\) consists of arrows \(Fj\to c\), and a colimit of \(F\) is a universal cone under \(F\). 

Despite referring to a limit as \emph{the} limit, limits are not unique on the nose.
They are unique up to isomorphism, and the isomorphism is unique among arrows which are compatible with the cones in the above way.

\begin{example}
    If \(J\) is discrete, that is, the only arrows are identities, then the limit of a diagram \(F\colon J\to \Cc\) is a \emph{product}, denoted \(\prod_{j\in J}Fj\), or \(Fj_1\times\cdots\times Fj_n\) if \(J\) is finite. 
    In many common categories (including \(\Set,\Top,\Grp,\Ring,\Vect_k\)) these coincide with the classical product.
\end{example}

\begin{example}\label{ex-equalisers}
    An \emph{equaliser} is a limit of a diagram where the shape category consists of two objects and two non-trivial parallel arrows from one object to another:
    \[\begin{tikzcd}
        \bullet & \bullet
        \arrow[shift right, from=1-1, to=1-2]
        \arrow[shift left, from=1-1, to=1-2]
    \end{tikzcd}\]
    If \(\varphi,\psi\colon G\to H\) are group homomorphisms, the \emph{equaliser} of \(\varphi\) and \(\psi\) is the limit of the diagram consisting of these two homomorphisms (its shape category consists of two objects and two non-trivial parallel arrows from one object to another).
    The equaliser determines the subgroup \(\{g\in G: \varphi(g)=\psi(g)\}\), and every subgroup arises this way.
    The kernel of a group homomorphism \(\varphi\) coincides with the equaliser between \(\varphi\) and the trivial homomorphism.
\end{example}

\begin{example}\label{ex-pullback}
    A \emph{pullback} is a limit of a diagram of shape
    \[\begin{tikzcd}
        & \bullet \\
        \bullet & \bullet
        \arrow[from=1-2, to=2-2]
        \arrow[from=2-1, to=2-2]
    \end{tikzcd}\]
    resulting in a limit
    \[\begin{tikzcd}
	{c\baltimes{f}{g}d} & d \\
        c & e
        \arrow[from=1-1, to=1-2]
        \arrow[from=1-1, to=2-1]
        \arrow["\lrcorner"{anchor=center, pos=0.125}, draw=none, from=1-1, to=2-2]
        \arrow["f", from=1-2, to=2-2]
        \arrow["g"', from=2-1, to=2-2]
    \end{tikzcd}\]
    In \(\Top\) the pullback is the usual pullback or \emph{fibre-product}: for \(f\colon X\to Z\) and \(g\colon Y\to Z\) the pullback \(X\baltimes{f}{g}Y\) is the subspace \(\{(x,y)\in X\times Y: f(x)=g(y)\}\) with the subspace topology coming from the product topology.
    The limit cone is given by the coordinate projections.
\end{example}

A category \(\Cc\) is \emph{concrete} (or rather, \emph{concretisable}) if there is a faithful functor \(U\colon \Cc\to\Set\).
We refer to \(U\) as a forgetful functor.
The underlying sets of the limits in Examples~\ref{ex-equalisers} and \ref{ex-pullback} in \(\Grp\) and \(\Top\) coincide with the sets you would get by taking the limits in \(\Set\) if you `forgot' the group or topological structures.
This is always the case: the usual forgetful functors \(\Grp,\Top\to\Set\) send limits to limits.

\begin{definition}
    Let \(G\colon\Cc\to\Dd\) be a functor, and consider a class \(\Cc\) of diagrams in \(\Cc\).
    The functor \(G\) 
    \begin{enumerate}
        \item \emph{preserves limits} in \(\Cc\), if for any limit cone \(\lambda\) over a diagram \(F\) in the class \(\Cc\), the image cone \(G\lambda\) is a limit over \(GF\),
        \item \emph{reflects limits} in \(\Cc\), if whenever the image \(G\lambda\) of a cone \(\lambda\) over a diagram \(F\) in \(\Cc\) is a limit cone over \(GF\), then \(\lambda\) is a limit cone over \(F\),
        \item \emph{creates limits} in \(\Cc\), if whenever \(GF\) has a limit cone \(\mu\) in \(\Dd\) for \(F\) in \(\Cc\), there is a cone \(\lambda\) over \(F\) with \(G\lambda=\mu\), and further, \(G\) reflects limits.
    \end{enumerate}
    We may also refer to \(G\) preserving, reflecting, or creating colimits in a similar fashion.
\end{definition}

\begin{example}
    The forgetful functors from \(\Top,\Grp,\Ring,\Vect\) to \(\Set\) all preserve limits.
    All except \(\Top\to\Set\) also create limits.
\end{example}

A category is \emph{complete} if it admits limits of all small diagrams.
Importantly, \(\Set\) is complete.
The limit of a small diagram \(F\colon J\to\Set\) is the set \(\Cone(\{\ast\},F)\) of cones over \(F\) where the apex is a singleton.
The \(j\)-th leg of the universal cone sends a cone \(\mu\colon\{\ast\}\to Fj\) to the element \(\mu_j(\ast)\). 
Given another cone \(\nu=(\nu_j\colon X\to Fj)_{j\in J}\), the universal map \(\tilde\nu\colon X\to \Cone(\{\ast\},F)\) sends \(x\in X\) to the cone \(\tilde\nu(x)_j\colon \ast\mapsto\nu_j(x)\).

Since \(\Set\) is complete, any category admitting a functor to \(\Set\) which creates small limits is also complete.

\subsection{Equivalences of categories and adjunctions}

\begin{definition}
    Let \(\Cc\) and \(\Dd\) be categories.
    An \emph{equivalence} between \(\Cc\) and \(\Dd\) consists of functors \(F\colon\Cc\to\Dd\) and \(G\colon \Dd\to\Cc\) such that there are natural isomorphisms \(GF\cong 1_\Cc\) and \(FG\cong 1_\Dd\).
    If there exists an equivalence between \(\Cc\) and \(\Dd\), we say \(\Cc\) and \(\Dd\) are \emph{equivalent}, denoted \(\Cc\simeq\Dd\).
    Assuming the axiom of choice, a functor \(F\colon\Cc\to\Dd\) gives rise to an equivalence if and only if
    \begin{enumerate}
        \item \(F\) is \emph{full} and \emph{faithful}: the map \(\colon\Cc(c,c')\to\Dd(Fc,Fc')\) induced by \(F\) is a bijection for all \(c,c'\in\Cc^0\),
        \item \(F\) is \emph{essentially surjective}: for every object \(d\in\Dd^0\) there is an object \(c\in\Cc^0\) such that \(Fc\cong d\).
    \end{enumerate}
\end{definition}

Categories can be equivalent without being isomorphic; the compositions \(FG\) and \(GF\) above may fail to be equal to the identity functors.
A further weakening of the notion of two functors being mutually inverse which remains very interesting is the following.

\begin{definition}
    Let \(F\colon \Cc\to\Dd\) and \(G\colon \Dd\to\Cc\) be functors.
    An \emph{adjunction} between \(F\) and \(G\) is a family of bijections
    \[\Dd(Fc,d)\xrightarrow\sim \Cc(c,Gd)\]
    which are natural in both \(\Cc\) and \(\Dd\).
\end{definition}

If such an adjunction exists, we say \(F\) is \emph{left-adjoint} to \(G\), or that \(F\) is the \emph{left adjoint} of \(G\).
We also say that \(G\) is \emph{right-adjoint} to \(F\), or that \(G\) is the \emph{right adjoint} to \(F\).
Left and right adjoints are unique up to natural isomorphism, justifying our use of the definite article.
We denote this with a turnstyle \(F\dashv G\).

One commonly occurring flavour of adjunctions consists of the ``free-forgetful'' pairs.
The following are some examples of these:
\begin{itemize}
    \item The free group functor \(\Set\to\Grp\) left-adjoint to the forgetful functor \(\Grp\to\Set\),
    \item The free vector space functor \(\Set\to\Vect\) is left-adjoint to the forgetful functor \(\Vect\to\Set\),
    \item Equipping a set with the discrete topology is left-adjoint to the forgetful functor \(\Top\to\Set\).
\end{itemize}

\begin{theorem}[{\cite[Theorem~V.5.1]{MacLane_WorkingMathematician}}]\label{thm-RAPL}
    Right adjoints preserve limits.
    Dually, left adjoints preserve colimits.
\end{theorem}

An adjunction \(F\dashv G\) yields natural transformations \(\eta\colon 1_\Cc\Rightarrow GF\) and \(\epsilon\colon FG\Rightarrow 1_\Dd\).
The component \(\eta_c\colon c\to GFc\) is the image of \(1_c\in\Cc(c,c)\) under the composition
\[\Cc(c,c)\xrightarrow{F}\Dd(Fc,Fc)\xrightarrow\sim \Cc(c,GFc),\]
and similarly the component \(\epsilon_d\colon FGd\to d\) is the image of \(1_d\in\Dd(d,d)\) under the composition
\[\Dd(d,d)\xrightarrow{G}\Cc(Gd,Gd)\xrightarrow\sim\Dd(FGd,d).\]
The natural transformations \(\eta\) and \(\epsilon\) are the \emph{unit} and \emph{counit} of the adjunction.
The triangle identities state that the compositions
\[\begin{tikzcd}
    F \arrow[r, Rightarrow, "F\eta"] \arrow[dr, Rightarrow, equal] & FGF \arrow[d, Rightarrow, "\epsilon F"] \\
    & F
\end{tikzcd}
\quad\text{and}\quad
\begin{tikzcd}
    G \arrow[r, Rightarrow, "\eta G"] \arrow[dr, Rightarrow, equal] & GFG \arrow[d, Rightarrow, "G\epsilon"] \\
    & G
\end{tikzcd}\]
are the identity natural transformations, where e.g. \(F\eta\colon F\to F(GF)\) is the natural transformation with components \((F\eta)_c=F\eta_c\) and \(\epsilon F\colon FGF\Rightarrow F\) is the natural transformation with components \((\epsilon F)_c=\epsilon_{Fc}\).
In fact, two functors form an adjunction exactly when there are natural transformations \(\eta\) and \(\epsilon\) satisfying the triangle identities.

Adjunctions restrict to equivalences on categories.
A \emph{subcategory} of a category \(\Cc\) consists of a subclass of objects and a subclass of morphisms which form a category with the same composition as \(\Cc\).
A subcategory \(\Dd\subseteq\Cc\) is \emph{full} if \(\Cc(c,d)=\Dd(c,d)\) for all objects \(c,d\in\Dd^0\subseteq\Cc^0\).

\begin{lemma}
    Let \(F\colon\Cc\to\Dd\) and \(G\colon\Dd\to\Cc\) form an adjunction \(F\dashv G\).
    The \emph{fixed points} in \(\Cc\) (resp. \(\Dd\)) are the objects 
    \(c\in\Cc^0\) (resp. \(d\in\Dd^0\)) for which the component \(\eta_c\colon c\to GFc\)  of the unit (resp. \(\epsilon_d\colon FGd\to d\) of the counit) is an isomorphism.
    The functors \(F\) and \(G\) restrict to an equivalence of categories between the full subcategories of \(\Cc\) and \(\Dd\) spanned by fixed points.
\end{lemma}

%% file: appendices/stoneDuality.tex
The topology \(\Oo X\) on a topological space \(X\) carries a canonical lattice structure for the partial ordering by inclusion; joins are given by unions and meets are given by intersections.
A continuous map between topological spaces also induces a lattice homomorphism between the topologies via preimaging: for a continuous function \(f\colon X\to Y\), the map \(\Oo Y\to\Oo X\) sending an open subset \(U\in\Oo Y\) to \(f^{-1}(U)\) is a lattice homomorphism. 
The lattice structure on a topology enjoys some further perks: it is closed under arbitrary unions (not just finite ones) and so the lattice has all suprema, and De Morgan's laws imply that meets distribute over arbitrary joins.

Frames are an abstract axiomatisation of the lattice structure and properties shared by all topologies.
Frames and their morphisms form a dual picture to topological spaces and continuous functions, and \emph{Stone duality} describes an adjunction between the two categories which restricts to an equivalence between sober topological spaces and `spatial' frames (i.e. frames coming from topological spaces).

We refer the reader to \cite{PicadoPultr_FramesLocales} for a comprehensive reference on the theory of frames and locales.

\begin{definition}
    A \emph{frame} \(F\) is a partially ordered set with finite infima and all suprema, in which infima distribute over suprema.
    We denote the infimum or \emph{meet} of two elements \(e,f\in F\) by \(e\wedge f\), and the supremum or \emph{join} of a family \((e_\alpha)_\alpha\subseteq F\) by \(\bigvee_\alpha e_\alpha\).
    Distributivity is then expressed as 
    \[e\wedge\bigvee_\alpha f_\alpha=\bigvee_\alpha (e\wedge f_\alpha)\]
    for all \(e,f_\alpha\in F\).

    Given frames \(F\) and \(L\), a \emph{frame homomorphism} \(\varphi\colon F\to L\) is a function which preserves finite infima and arbitrary suprema. 
    Frames together with frame homomorphisms form the category \(\Frame\).
\end{definition}

We quickly note that frame homomorphisms are necessarily order-preserving since the order in a lattice can be expressed in terms of meets or joins.
Explicitly, for \(a,b\in F\) we have \(a\leq b\) if and only if \(a\vee b=b\), if and only if \(a\wedge b=a\).

\begin{remark}\label{rem-meetPresNecc}
    The requirement that a frame homomorphism preserves finite infima is necessary for our purposes and it does not follow automatically from preservation of suprema.
    An order-preserving map \(f\colon A\to B\) between lattices will always satisfy \(f(a\wedge b)\leq f(a)\wedge f(b)\), but we do not always gain equality, even if \(f\) is a map between frames which preserves suprema.
    We construct an explicit counterexample.
    Consider the frame \(\Omega=\{0<1\}\) and equip \(\Omega^2\) with the product order.
    Thinking of \(\Omega\) as a set of truth values, the logical `or' function \(+\colon\Omega^2\to\Omega\), \(+(a,b)=\max\{a,b\}\) is order preserving and preserves all suprema, since the only families of values \((a_\alpha,b_\alpha)_\alpha\) whose join is sent to zero under \(+\) are families consisting only of zeroes.
    But \(+\) does not preserve binary meets, since \(+((1,0)\wedge(0,1))=+(0,0)=0\), but \(+(1,0)\wedge +(0,1)=1\wedge 1=1\).
\end{remark}

Note that a frame \(F\) has both an absolute minimum (often denoted by \(0\)) and an absolute maximum (denoted by \(1\)) since frames have finite suprema and infima, and the empty join and meet are each respectively a minimum and maximum for the frame. 
These may also be called the \emph{bottom} or \emph{zero} (respectively \emph{top} or \emph{1}) of the frame.
It follows by definition that frame homomorphisms preserve the minimum and maximum of a frame.

As previously introduced, the first motivating examples of frames we encounter are topologies.
A topological space is, by definition, a set \(X\) paired with a collection \(\Oo X\) of `open' subsets which is closed under finite intersections and arbitrary unions. 
When equipped with the partial ordering of inclusion, the axioms of a topological space immediately show that \(\Oo X\) is a frame; the meet is given by intersection and the join is given by union.
A continuous function \(f\colon X\to Y\) between topological spaces is exactly a function for which the induced preimaging map \(f^{-1}\) forms a function \(\Oo Y \to\Oo X\).
Preimaging preserves all unions and intersections, so \(f^{-1}\) is a frame homomorphism, and we quickly see that \(\Oo\) froms a contravariant functor from \(\Top\) to \(\Frame\).

Somewhat surprisingly, under quite mild hypotheses it is possible to recover the underlying spaces \(X\) and \(Y\), as well as the function \(f\) from the frame structures of \(\Oo X\) and \(\Oo Y \) together with the homomorphism \(f^{-1}\).
This reconstruction is the first obvious raison d'\^{e}tre for the study of sober spaces and Stone duality.
Sober spaces, as we shall see, are precisely those for which this reconstruction functions.

\begin{definition}{{\cite[1.1]{PicadoPultr_FramesLocales}}}\label{defn-irredSob}
    A closed subset \(C\) of a topological space \(X\) is \emph{irreducible} if it is non-empty and for any closed sets \(C_1,C_2\subseteq X\) with \(C_1\cup C_2=C\), we have \(C_1=C\) or \(C_2=C\).
    We say \(X\) is \emph{sober} if each closed irreducible subset is the closure of a unique singleton.

    We denote by \(\Top_\Sob\) the full subcategory of \(\Top\) spanned by sober spaces.
\end{definition}

One quickly notes that the closure of any singleton is necessarily an irreducible subset in any topological space, and in any \(T_0\) distinct singletons will yield distinct closures.
It follows that sober spaces are always \(T_0\), which is to be expected given the earlier claim that a sober space can be recovered from the lattice structure of its topology, and non-\(T_0\) spaces contain pairs of points which cannot possibly be distinguished by the topology.

The difference between sober and \(T_0\) is that singletons are the only way in which closed irreducible subsets may arise. 

Sober spaces are abundant: firstly one notes that any Hausdorff space is sober.
Sobriety also passes `up' local homeomorphisms in the following sense, which will be relevant later in the discussion of \'etale groupoids.

\begin{lemma}
    All Hausdorff spaces are sober, and all sober spaces are \(T_0\).
    Sober and \(T_1\) are inequivalent.
    \begin{proof}
        The first claim is exactly the above discussion.
        We provide two examples: a sober space which is not \(T_1\) and a \(T_1\) space which is not sober.

        Firstly, a space that is sober but not \(T_1\).
        Consider the Sierpinski space \(X=\{0,1\}\) with the topology \(\Oo X=\{\emptyset,\{0\},X\}\).
        The closed irreducible subsets are exactly \(X\) and \(\{1\}\).
        We note that \(X\) is the closure of \(\{0\}\) and \(\{1\}\) is the closure of itself, so \(X\) is sober. 
        The singleton \(\{0\}\) is however not closed, so \(X\) fails to be \(T_1\).

        Next, a \(T_1\) space which is not sober.
        Consider the natural numbers \(\NN\) equipped with the cofinite topology i.e. the closed sets are exactly the finite subsets of \(\NN\) and \(\NN\) itself.
        Then each singleton is closed, so this is a \(T_1\) space.
        The closed set \(\NN\) is itself irreducible, since if the finite union of closed sets is either finite or one of the closed sets is already \(\NN\).
        But \(\NN\) is not the closure of any singleton, so this space fails to be sober.
    \end{proof}
\end{lemma}

The forgetful functor \(\Top\to\Set\) sending a topological space to its underlying set is represented by the singleton space \(\{\ast\}\); the usual choice of natural bijections send a point \(x\in X\) to the 'selecting' function \(f_x\colon \{\ast\}\to X\), \(f_x(\ast)=x\).
Dualising to preimages, the frame homomorphisms \(f_x^{-1}\colon\Oo X\to\Oo \{\ast\} \) send an open set \(U\in\Oo X\) to \(\{\ast\}\) if \(x\in U\), and to \(\emptyset\) otherwise. 
Moreover, if \(X\) is \(T_0\) then the induced preimaging maps distinguish points. 
We shall see that sober spaces are those for which every such `character' on \(\Oo X\) corresponds exactly to a point in \(X\).

\subsection{Stone duality}

Using the universal property of the singleton topological space as inspiration, we reconstruct an underlying set and topology from a given frame \(F\). 
This space can always be constructed, if \(F\cong\Oo X\) for a sober space \(X\) then the reconstructed space is homeomorphic to \(X\).
We denote by \(\Omega=\{0<1\}\) the unique totally ordered two-element frame.

\begin{definition}
    A \emph{character} on a frame \(F\) is a frame homomorphism \(\chi\colon F\to\Omega\).
    Denote the set of characters by \(\widehat F\).
\end{definition}

There is a canonical topology to place on \(\widehat F\).

\begin{lemma}\label{lem-frameSpecTopo}
    The sets \(\Uu_e:=\{\chi\in\widehat F:\chi(e)=1\}\) ranging over \(e\in F\) define a topology on \(\widehat F\), and the map \(\eta_F\colon F\to\Oo \widehat F \) given by \(\eta_F(e)=\Uu_e\) is a surjective frame homomorphism.
    \begin{proof}
        The image of a frame homomorphism is itself clearly a frame so it suffices to show that \(\eta_F\) preserves joins and finite meets.
        First note that \(\Uu_1=\widehat F\) and \(\Uu_0=\emptyset\) since every character preserves absolute minima and maxima.
        For any family \((e_\alpha)_\alpha\subseteq F\) write \(e:=\bigvee_\alpha e_\alpha\).
        Since they preserve all joins, any character \(\chi\in\widehat F\) satisfies \(\chi(e)=1\) if and only if \(\chi(e_\alpha)=1\) for some \(\alpha\), so \(\Uu_e =\bigcup_\alpha\Uu_{e_\alpha}\).
        For \(e,f\in\widehat F\) we have \(1=\chi(e\wedge f)=\chi(e)\wedge\chi(f)\) if and only if \(\chi(e)=\chi(f)=1\), so \(\Uu_e\cap\Uu_f=\Uu_{e\wedge f}\).
        Hence \(\eta_F\) preserves all joins and finite meets, whereby it is a frame homomorphism.
        It is surjective by construction.
    \end{proof}
\end{lemma}

The space \(\widehat F\) is called the \emph{spectrum} of \(F\), and we always consider it equipped with the topology from Lemma~\ref{lem-frameSpecTopo}.
Since the composition of frame homomorphisms is again a frame homomorphism, any frame homomorphism \(\varphi\colon F\to L\) induces a map \(\varphi^*\colon \widehat L\to\widehat F\) by precomposition (remembering that elements of \(\widehat L\) and \(\widehat F\) are frame homomorphisms into \(\Omega\)).
This function is continuous.

\begin{lemma}\label{lem-precompContNat}
    Let \(\varphi\colon F\to L\) be a frame homomorphism.
    The preimage of the open subset \(\Uu_e\subseteq\widehat F\) under the function \(\varphi^*\colon \widehat L\to\widehat F\) is exactly \(\eta_L(\varphi(e))=\Uu_{\varphi(e)}\).
    In particular, \(\varphi^*\) is continuous.
    \begin{proof}
        By Lemma~\ref{lem-frameSpecTopo}, a generic open subset of \(\widehat F\) has the form \(\eta_F(e)=\Uu_e\) for some \(e\in F\).
        We have \(\chi\in\Uu_{\varphi(e)}\) exactly when \(1=\chi(\varphi(e))=\varphi^*\chi(e)\), exactly when \(\varphi^*\chi\in\Uu_e\).
        This yields \((\varphi^*)^{-1}(\Uu_e)=\Uu_{\varphi(e)}\).
    \end{proof}
\end{lemma}

We clearly have \((\varphi\psi)^*=\psi^*\varphi^*\) for composable frame homomorphisms \(\varphi\) and \(\psi\), as well as \(1_F^*=1_{\widehat F}\) (where \(1_c\) is the identity morphism for the object \(c\) in a category), so we get a functor \(\sigma\colon\Frame\to\Top^\op\) sending a frame \(F\) to its spectrum \(\sigma F:=\widehat F\) and a frame homomorphism \(\varphi\) to precomposition with it \(\sigma\varphi :=\varphi^*\).
Moreover we have shown:

\begin{corollary}\label{cor-etaFrameNatural}
    The frame homomorphisms \(\eta_F\colon F\to\Oo \sigma F \) assemble into a natural transformation \(1_\Frame\Rightarrow \Oo \sigma\).
\end{corollary}

The spectrum \(\sigma F\) of a frame \(F\) is always sober.
\begin{proposition}[cf. {\cite[II~1.7]{Johnstone_StoneSpaces}}]\label{prop-specOfFrmSob}
    Let \(F\) be a frame.
    Then \(\sigma F\) is sober.
    \begin{proof}
        First note that characters \(\sigma F\) are distinct only if they differ on some element \(f\in F\), and hence the open set \(\Uu_f\subseteq \sigma F\) separates them.
        Thus \(\sigma F\) is \(T_0\) and thus two singletons have the same closure only if they are equal.

        Fix an irreducible closed subset \(C\subseteq \sigma F\) and let \(U:=\sigma F\setminus C\subseteq \sigma F\) be its complement.
        Since \(U\) is open, Lemma~\ref{lem-frameSpecTopo} yields some \(e\in F\) such that \(\Uu_e=U\), and note that \(e\neq 1\) since \(\Uu_e\neq \sigma F\).
        Moreover, by replacing \(e\) with the join \(\bigvee\{f\in F: \Uu_f\subseteq U\}\), we may assume that \(e\) is the maximal such element of \(F\) with \(\Uu_e=\sigma F\setminus C\), and that \(\Uu_f\subseteq\Uu_e\) implies \(f\leq e\) for any \(f\in F\).
        For any \(f_1,f_2\in F\) with \(f_1\wedge f_2\leq e\) we then have \(\Uu_{f_1\wedge f_2}\subseteq\Uu_e\), whereby \(\Uu_{f_1}\subseteq\Uu_e\) or \(\Uu_{f_2}\subseteq \Uu_e\) since the complement of \(\Uu_e\) is irreducible.
        But since \(e\) is chosen as the above join, we necessarily have \(f_1\leq e\) or \(f_2\leq e\).
        
        Define \(\chi\colon F\to\Omega\) by
        \[\chi(f)=\begin{cases}
            1,&f\nleq e,\\
            0,&f\leq e.
        \end{cases}\]
        We claim that \(\chi\) is a character and \(C=\overline{\{\chi\}}\).
        We have \(\chi(0)=0\) as \(0\leq e\), and \(\chi(1)=1\) since \(1\nleq e\).
        For a family \((f_\alpha)_\alpha\in F\), the join \(\bigvee_\alpha f_\alpha\) is mapped to 0 if and only if \(f_\alpha\leq e\) for each \(\alpha\) i.e. \(\chi(f_\alpha)=0\) for each \(\alpha\).
        Hence \(\chi(f)=0\) if and only if \(\chi(f_\alpha)=0\) for all \(\alpha\), and so \(\chi\) preserves joins.
        Lastly, \(\chi\) preserves meets since \(\chi(f_1\wedge f_2)=0\) if and only if \(f_1\wedge f_2\leq e\).
        The choice of \(e\) being maximal together with the previous discussion then yields \(f_1\leq e\) or \(f_2\leq e\), and hence \(\chi(f_1)\wedge\chi(f_2)=0\).
        The reverse implication is clear.
        Thus \(\chi\) is a character, and by construction \(\chi(e)=0\) so \(\chi\notin\Uu_e\).
        
        We have \(\varsigma\in C\) if and only if \(\varsigma(e)=0\).
        Let \(\Uu_f\) be an open neighbourhood of \(\varsigma\).
        Then \(\varsigma(e)=0<1=\varsigma(f)\) and so \(f\nleq e\), since \(\varsigma\) preserves order.
        Hence \(\chi(f)=1\) and so \(\chi\in\Uu_f\), whereby \(\{\chi\}\) is dense in \(C\).
        This complete the proof.
    \end{proof}
\end{proposition}

Considering the other direction: given a topological space \(X\) there is a natural map \(\kappa_X\colon X\to S\Oo X\) sending an point \(x\in X\) to the character \(\kappa_X(x):=\chi_x\in S\Oo X\) which detects it:
\[\chi_x(U)=\begin{cases}
    1,& x\in U,\\
    0,& x\notin U.
\end{cases}\]
Sober spaces are those for which this map is a homeomorphism; allowing us to recover \(X\) from the frame structure of its topology.

\begin{theorem}[cf. {\cite[6.2~Proposition]{PicadoPultr_FramesLocales}}]\label{thm-soberIffSpecHomeo}
    For any topological space \(X\), the map \(\kappa_X\) is continuous.
    It is injective if and only if \(X\) is \(T_0\), and it is a homeomorphism if and only if \(X\) is sober.
    \begin{proof}
        By Lemma~\ref{lem-frameSpecTopo}, a generic open subset of \(\sigma \Oo X\) has the form \(\Uu_U=\{\chi\in \sigma \Oo X:\chi(U)=1\}\) for an open subset \(U\in\Oo X\). 
        The preimage \(\kappa_X^{-1}(\Uu_U)\) consists of points \(x\) for which \(\chi_x(U)=1\), which are exaclty the points in \(U\). 
        Thus \(\kappa_X\) is continuous.
        
        If \(X\) is \(T_0\) then for any two (distinct) points in \(X\) there is an open set which contains one and not the other.
        The resulting characters then take different values on this open subset, so \(\kappa_X\) is injective.
        Conversely, if \(\kappa_X\) is injective then for any two distinct points \(x,y\in X\) there is an open subset \(U\in\Oo X\) with \(\chi_x(U)\neq\chi_y(U)\).
        This says precisely that \(U\) contains one of the two points \(x,y\) but not the other, which is to say that \(X\) is \(T_0\).

        The space \(\sigma \Oo X\) is sober by Proposition~\ref{prop-specOfFrmSob}, so if \(\kappa_X\) is a homeomorphism then \(X\) is also sober.

        Now suppose that \(X\) is sober.
        Any character \(\chi\in \sigma \Oo X\) preserves all joins, so there is a largest open subset \(U\Oo X\) with \(\chi(U)=0\) (namely the join of all open subsets which \(\chi\) sends to zero).
        For any open subsets \(U_1,U_2\in\Oo X\) we have \(U_1\cap U_2\subseteq U\) if and only if \(\chi(U_1)\wedge\chi(U_2)=\chi(U_1\wedge U_2)=0\), whereby one of \(U_1\) or \(U_2\) is already contained in \(U\).
        This exactly says that \(C=X\setminus U\) is a closed irreducible subset, and so is the closure of a unique singleton \(x\in X\) by sobriety of \(X\).
        For any open subset \(V\in\Oo X\) we then have \(\chi(V)=1\) if and only if \(V\) is not contained in \(U\).
        But \(V\) then has non-empty intersection with \(C=\overline{\{x\}}\), whereby \(x\in V\) since \(V\) is open.
        Thus \(\chi(V)=1\) implies \(\chi_x(V)=1\), equivalently \(x\in V\), and obviously we have \(\chi(V)=1\) if \(x\in V\), so we see \(\chi=\chi_x\).
        Thus \(\kappa_X\) is surjective.

        It remains to show that \(\kappa_X\) is open.
        Since \(\kappa_X\) is surjective, we know that \(\sigma \Oo X=\{\chi_x:x\in X\}\).
        The image of an open set \(U\in\Oo X\) under \(\kappa_X\) is then the set \(\{\chi_x: x\in U\}=\{\chi_x:x\in X, \chi_x(U)=1\}=\Uu_U\), which is open in \(\sigma \Oo X\) as required.
    \end{proof}
\end{theorem}

Thus we see that sober topological spaces are exactly those which may be systematically recovered by their frames of open sets in this way. 
The composition of \(\sigma \) and \(\Oo\) yields a functor \(\Top\to\Top\), and the maps \(\kappa_X\colon X\to \sigma \Oo X\) are kind enough to be natural in \(X\).

\begin{lemma}
    The functions \(\kappa_X\colon X\to \sigma \Oo X\) assemble into a natural transformation \(1_\Top\Rightarrow \sigma \Oo\).
    \begin{proof}
        Fix a continuous map \(f\colon X\to Y\).
        The induced map \(\sigma  \Oo f\colon \sigma \Oo X\to \sigma \Oo Y\) is \((f^{-1})^*\); precomposition with the preimaging function associated to \(f\).
        This sends a character \(\chi\in \sigma \Oo X\) to \((f^{-1})^*\chi=\chi\circ f^{-1}\).
        We compute for \(x\in X\) and \(U\in\Oo Y \):
        \[[\kappa_Y\circ f(x)](U)=\chi_{f(x)}(U)=\chi_x(f^{-1}(U))=\kappa_X(x)\circ f^{-1}(U)=[(f^{-1})^*\kappa_X(x)](U),\]
        yielding naturality.
    \end{proof}
\end{lemma}

While shall not consider such objects, there do exist frames which are not isomorphic to any topology.
Such `non-spatial' frames will not be treated in detail in this article, but their existence warrants the exercise of some caution when attempting to wield Stone duality.
We refer the curious reader to other resources (e.g. \cite{Johnstone_Elephant1}).

\begin{definition}
    A frame is \emph{spatial} if it is isomorphic to a topology.
    The full subcategory of \(\Frame\) spanned by spatial frames is denoted \(\Frame_\Spat\).
\end{definition}

A priori the spatiality of a frame \(F\) can be realised by any topology, but by \cite[5.1~Proposition]{PicadoPultr_FramesLocales} a frame is spatial exactly when the homomorphism \(\eta_F\colon F\to\Oo \sigma F \) is an isomorphism. 
Combining Proposition~\ref{prop-specOfFrmSob} and Theorem~\ref{thm-soberIffSpecHomeo} we see that \(F\) is spatial if and only if it is isomorphic to the topology of its spectrum, and this is the only sober space with this property.

The version of Stone duality we shall present and use is given in the context of sober topological spaces.
The functors \(\sigma\) and \(\Oo\) form an adjunction which restricts to an equivalence between the categories of sober topological spaces and spatial frames.
This is a standard result, and we provide a proof for convenience.

\begin{theorem}\label{thm-appSoberSpatAdj}
    The functors \(\sigma \colon\Frame\to\Top^\op\) and \(\Oo\colon\Top^\op\to\Frame\) are an adjunction pair \(\sigma \dashv \Oo\) with unit \(\eta\colon 1_\Frame\Rightarrow \Oo \sigma \) and counit \(\kappa_X^\op\colon \sigma \Oo\Rightarrow 1_{\Top^\op}\).
    The fixed points of the adjunction are spatial frames and sober topological spaces, yielding an equivalence \(\Frame_{\Spat}\simeq\Top^\op_\Sob\).
    \begin{proof}
        The triangle identities to verify are \((\kappa^\op \sigma )(\sigma \eta)=1_\sigma \) and \((\Oo\kappa^\op)(\eta\Oo)=1_\Oo\).

        For a frame \(F\in\Frame\), the component of \((\kappa^\op \sigma )(\sigma \eta)\) at \(F\) is given by \(\kappa^\op_{\sigma F}\sigma \eta_{\sigma F}\), which is a arrow \(\sigma F\to \sigma F\) in \(\Top^\op\) dual to the function \(\eta_{\sigma F}^*\circ\kappa_{\sigma F}\colon \sigma F\to \sigma F\).
        For a character \(\chi\in \sigma F\), the resulting image of \(\chi\) under this map sends a frame element \(e\in F\) to
        \begin{align*}
            [\eta_{\sigma F}^*\circ\kappa_{\sigma F}(\chi)](e)&=[\kappa_{\sigma F}(\chi)](\eta_{\sigma F}(e))\\
            &=[\kappa_{\sigma F}(\chi)](\Uu_e)\\
            &=\begin{cases}
                1,&\chi\in\Uu_e\\
                0,&\chi\notin\Uu_e,
            \end{cases}\\
            &=\chi(e).
        \end{align*}
        Since the choice of \(F\), \(\chi\), and \(e\) were arbitrary, we have \((\kappa^\op\sigma )(\sigma \eta)=1_\sigma \).

        For the second triangle identity, for a topological space \(X\), the component of \((\Oo\kappa^\op)(\eta\Oo)\) at \(X\) is the frame homomorphism \(\kappa_X^{-1}\circ\eta_{\Oo X}\colon \Oo X\to \Oo X\).
        The first map \(\eta_{\Oo X}\) sends an open set \(U\in\Oo X\) to the set \(\Uu_U:=\{\chi\in \sigma \Oo X: \chi(U)=1\}\), and the second map sends \(\Uu_U\) to the preimage \(\kappa_X^{-1}(\Uu_U)\) which is equal to \(U\) by prior reasoning.

        Lastly, we combine Proposition~\ref{prop-specOfFrmSob} and Theorem~\ref{thm-soberIffSpecHomeo} to see that \(\kappa_{\sigma F}\) is a homeomorphism for all frames \(F\), so \(\kappa \sigma \) is a natural isomorphism. 
        It is clear that spatial frames and sober spaces are the respective fixed points of these functor, so we gain the desired equivalence of categories.
    \end{proof}
\end{theorem}

By employing the equivalence of categories from Theorem~\ref{thm-soberSpatAdj}, we may move between the topological and frame-theoretic contexts as we see fit.
This also shows that the full subcategory of spatial frames is a \emph{reflective subcategory} of \(\Frame\), meaning that the inclusion functor admits a left adjoint.
Analogously the full (opposite) subcategory of sober topological spaces is a \emph{coreflective subcategory} of \(\Top^\op\), meaning that the inclusion functor admits a right adjoint.
Dualising shows that sober spaces are also a reflective subcategory of \(\Top\).